\documentclass[10pt]{article}
\usepackage{hyperref}
\usepackage{xcolor}
\usepackage{amssymb}
\usepackage{amsmath}
\usepackage{cancel}
\usepackage{multicol}
\input{liemacs10.sty} 
\usepackage{color}

\renewcommand{\down}{{\mathop{\downarrow}}}
\newcommand\be{{\bf{e}}}

\renewcommand\up{{\uparrow}}

\newcommand{\sH}{{\sf H}}
\newcommand{\sK}{{\sf K}}
\newcommand{\sN}{{\sf N}}
\newcommand{\sV}{{\tt V}}

\newcommand\RR{{\mathbb R}}

\newcommand\bC{{\mathbb C}}

\newcommand{\Mob}{{\rm\textsf{M\"ob}}}

\renewcommand{\bO}{\mathbb O}

\renewcommand{\phi}{\varphi}

\newcommand{\Stand}{\mathop{{\rm Stand}}\nolimits}

\renewcommand\mlabel{\label}

\begin{document}

\title{A geometric perspective on \\ Algebraic Quantum Field Theory}
\author{Vincenzo Morinelli\\Dipartimento di Matematica, Universit\`a di Roma ``Tor
Vergata''\\ 
{\tt morinell@mat.uniroma2.it}}

\date{}
\maketitle

\begin{center}\emph{Dedicated to Karl-Hermann Neeb on the occasion of his sixtieth birthday.}
\end{center}

\begin{abstract}
In this paper we give a streamlined overview of some of the recent constructions provided with K.-H. Neeb, G. \'Olafsson and collaborators for a new geometric approach to Algebraic Quantum Field Theory.
Motivations, fundamental concepts and some of the relevant results about the abstract structure of these models are here presented. 
\end{abstract}

\tableofcontents

\section{Introduction}

These notes aim to provide a structured guideline, building upon recent results developed in collaboration with K.-H. Neeb and G. Ólafsson, to explore a geometric perspective on Algebraic Quantum Field Theory. In this presentation, we will focus on some selected outcomes from the literature of recent years, illustrating them through specific cases. We hope this approach will assist readers in navigating the main findings and the research directions, while maintaining a balanced view of both the geometric and operator-algebraic perspectives on the related topics.

We start our discussion with some  motivations.  Algebraic  Quantum Field Theory (AQFT)  is an axiomatic framework for Quantum Field Theory where a model is  specified in the sense of the Haag-Kastler net by a map associating to any open region $\cO$ of the spacetime $M$ its von Neumann algebra $\cA(\cO)$ representing the set of local observables, acting on a fixed complex Hilbert  space $\mathcal H$ (the state space), satisfying fundamental quantum and relativistic assumptions as:
\begin{itemize}
\item Isotony: 
  $\cO_1 \subeq \cO_2$ implies $\cA(\cO_1) \subeq \cA(\cO_2)$ ;
\item Locality: 
  $\cO_1 \subeq \cO_2'$ implies $\cA(\cO_1) \subeq \cA(\cO_2)'$,
  where $\cO'$ is the ``causal complement'' of $\cO$, i.e., the maximal
  open subset that cannot be connected to $\cO$ by causal curves;
\item Covariance: 
  There is a Lie group $G$ acting on $M$ and a unitary representation
  $U \: G \to \U(\cH)$ such that 
  $U_g \cA(\cO) U_g^{-1} = \cA(g\cO)$ for $g \in G$;
\item Uniqueness and Invariance of vacuum: there exists a unique (up to a phase) unit vector $\Omega\in\cH$ such that $U(g)\Omega = \Omega$   for every $g \in G$. $\Omega$ is called the vacuum state;
\item Reeh--Schlieder property:   $\Omega\in \cH$  is cyclic for $\cA(\cO)$ if $\cO \not=\eset$
\end{itemize}
(cf.~\cite{Ha96}). Note that given an open region $\cO_1$, if there exists an open region $\cO_2\subset\cO_1'$ then, by locality and Reeh-Schlieder property, $\Omega$ is also separating for $\cA(O_1)$. 
In this framework there is a deep interplay between the geometric and the algebraic parts of the theory.

The Tomita-Takesaki theory states that given a von Neumann algebra $\cA\subset\cB(\cH)$ with a cyclic and separating vector  $\Omega\in\cH$, one can define the \textit{Tomita operator }$S_{\cA,\Omega}$ as the closure  of the anti-linear involution $$\cA\Omega \ni a\Omega\longmapsto a^*\Omega\in\cA\Omega.$$  The polar decomposition $S_{\cA,\Omega}=J_{\cA,\Omega}\Delta^{\frac12}_{\cA,\Omega}$ defines an antiunitary involution $J_{\cA,\Omega}$ called \textit{modular involution} and a self-adjoint positive operator $\Delta_{\cA,\Omega}$ called \textit{modular operator} satisfying $$J_{\cA,\Omega}\cA J_{\cA,\Omega}=\cA',\quad \Delta_{\cA,\Omega}^{it}\cA\Delta_{\cA,\Omega}^{-it}=\cA,\quad J_{\cA,\Omega}\Delta_{\cA,\Omega} J_{\cA,\Omega}=\Delta_{\cA,\Omega}^{-1},$$ The one--parameter group $\Delta_{\cA,\Omega}^{it}$ is called the \textit{modular group}.  
It is relevant to remark that couples of an antiunitary operator $J$ and a self-adjoint positive operator $\Delta$ on a complex Hilbert space $\cH$ satisfying $J\Delta J=\Delta^{-1}$ are in 1-1 correspondence with standard subspaces $\sH\subset\cH$, namely real closed subspaces $\sH\subset\cH$ such that $\overline{\sH+i\sH}=\cH$ and $\sH\cap i\sH=\{0\}$, through the condition \begin{equation}\label{eq:HS}\sH=\ker(1-S)\end{equation} where $S=J\Delta^{\frac12}$ (cf. Sect.~\ref{sect:stsub}).
Because of this relation with the Tomita operator, the language of standard subspaces has a key role in our analysis.

Now, consider a net of von Neumann algebras $\cA$ on  the $1+d$ dimensional Minkowksi space $\RR^{1,d}$. Let $\cA(W_R)$ be von Neumann algebra of observables associated to the Rindler wedge
\begin{equation}\label{eq:Rinwed}W_R:=\{x=(x_0,\ldots,x_d)\in \RR^{1,d}: |x_0|<x_1\},\end{equation}
The \textit{Bisognano-Wichmann property} states that the modular group of  $\cA(W_R)$ with respect to the vacuum state $\Omega$ unitarily implements the following one--parameter subgroup of the Poincaré group, known as Lorentz boosts:
\begin{equation}\label{eq:boost}
\Lambda_{W_R}(t)=\left(\begin{array}{ccc}
\cosh(t)&\sinh(t) &\textbf{0}\\
\sinh(t)&\cosh(t)&\textbf{0}\\
\textbf{0}&\textbf{0}&\textbf{1}
\end{array}\right).
\end{equation}
Specifically, this implementation is given by
$$\Delta_{\cA(W_R),\Omega}^{it}=U(\Lambda_{W_R}(-2\pi t)).$$ 
The modular conjugation also has a geometric action when the \textit{modular reflection property }is satisfied:
$$U(j_{W_R})=J_{\cA(W),\Omega}$$
where $j_{W_R}=\diag(-1,-1,\textbf{1}_{d-1})$ reflects $j_{W_R}W_R=W_R'$. 
It is evident that the Bisognano-Wichmann property together with the modular reflection property (see also the PCT Theorem, Section \ref{sect:ff}) determine a geometric action of the modular group and the modular conjugation of the local algebras associated with certain specific regions.
They have been verified for a large number of models  (see e.g.~\cite{BW75,BGL93,Mu01, Mo18, DM20}), ensuring a very strong relation between the geometry and the algebraic structures, which enriches the algebraic approach to Quantum Field Theory in many directions. For instance, in the recent years, operator algebraic methods have proven to be highly effective in studying energy conditions of AQFT models.  Here the modular Hamiltonian, namely  the logarithm of the  modular operator, associated with the local algebra of  wedge regions  plays a central role, see for instance \cite{MTW22, Lo20, CLRR22, CLR20, Ara76, LX18, LM23, LM24, Wit18}. By the Bisognano--Wichmann property it has a geometrical meaning and  a rigorous analysis from both the algebraic and the geometric perspectives is possible.

One can establish a 1-1 correspondence between the Poincar\'e transforms of the Rindler wedge called wedge regions  $W\subset\RR^{1,d}$  and one--parameter group of boosts $\Lambda_W$ obtained by Poincar\'e adjoint action on $\Lambda_{W_R}$. The algebraic canonical construction of the fundamental model called the free field provided by Brunetti–Guido–Longo (BGL) builds on  this wedge-boost correspondence, the Bisognano--Wichmann (BW) property and the modular reflection property, cf.~\cite{BGL02}.  Here, the standard subspace techniques play a prominent role. Given a (anti-)unitary representation $U$ of the Poincar\'e group (a particle in QFT) the localization properties of particle states (vectors of the Hilbert space $\cH$ supporting $U$)  on wedge regions can be characterized using nets of standard subspaces:
There is a map associating to each wedge region $W$ a standard subspace $\sH(W)\subset\cH$ of vector states in a canonical way by assuming the Bisognano-Wichmann property, modular reflection property, and taking into account the correspondence between the Tomita operators and standard subspaces given by \eqref{eq:HS}. The free net of von Neumann algebras is then constructed by a canonical procedure called second quantization.
 This approach is actually very general and this is already visible in the BGL paper. Indeed, the authors  apply this construction to many of  the free models  of physical interest  that have an analogue 1-1 correspondence between specific regions of spacetime  and one--parameter subgroups of the symmetry group:  In free (possibly conformal) theories on Minkowksi spacetime, on de Sitter space and on the chiral circle (see also \cite{GLW,Lo08}).

In this work we present a geometric perspective on these examples that strictly generalizes the framework given by the AQFT.
The paper starts giving more details on the wedge-boost generator correspondence in some of the AQFT spacetime (Sect.~\ref{sect:wre}).
We will explore how this picture can be generalized using a geometric approach, following \cite{MN21}. The Lie algebra elements generating  the one-parameter groups of boosts that identify wedge regions are Euler elements: Let $x$ be an element of $ \fg$ a finite dimensional real Lie algebra, $x$ is an {\it Euler element} if $\ad x$ is diagonalizable with 
$\Spec(\ad x) \subeq \{-1,0,1\} $
(see \cite{BN04, Kan00} for more details  on Euler elements and three graded Lie algebras). In Section \ref{sect:ee} we will see how this definition leads to a notion of  abstract wedges, called Euler wedges,  for a $\mathbb Z_2$-graded Lie group whose Lie algebra contains Euler elements without necessity to refer to a specific manifold. Simple real Lie algebras that support Euler elements, can be classified \cite[Thm.~3.10]{MN21}. 
One can study causal homogenous spaces to recover a manifold playing the role of the spacetime. In Section \ref{sect:chs} we present a suitable definition of wedge regions on causal homogeneous spaces and as an example we will see how the de Sitter spacetime can be recovered in this general picture (cf.~\cite{NO23,NO22,MNO23a,MNO23b}).

This allows us to exhibit a very general setting for nets of standard subspaces and von Neumann algebras having as an index set the Euler wedges. We will see in Section \ref{sect:geoss} that this construction is possible for a large family of Lie groups and provides several new models. 
An axiomatic framework is also provided for nets of standard subspaces on homogeneous spaces see \cite{FNO23, MNO23a, MNO23b,NOO21,NO21,NO22}. 
In Section \ref{sect:modEul} we discuss a recent construction of non-modular covariant net of standard subspaces, namely nets of standard subspaces with no geometric action of the wedge subspace modular groups  as in \cite[Sect.~3.1]{MN22}. In Section \ref{subsec:3.1} we introduce the Euler Element Theorem that ensures that it is possible to deduce  the keyrole of  Euler elements and Euler wedges from the Bisognano--Wichmann property and a regularity property. 
The latter property is strongly related with a localization property, as reported in  sections \ref{subsec:3.1} and \ref{sec:loc-reg} and in \cite[Sect.~4]{MN24}. Lastly we show that the von Neumann algebras appearing in this picture are of type III$_1$ for the Connes classification, see Section \ref{sect:vna} and \cite[Sect.~5]{MN24}. Many of these results are formulated for general von Neumann algebras and/or standard subspaces in some specific  relative position. We explain how they apply to nets of standard subspaces and of von Neumann algebras. This approach provides feedback for representation theory (see e.g.  constructions of nets of standard subspaces on homogeneous spaces in \cite{FNO23, NO21}) and for the algebraic approach to Quantum Field Theory without restrictions to second quantization models. This paper will mainly deal with abstract wedges, a more detailed overview on the properties of wedge regions in causal symmetric space is postponed.

\medskip

\nin{\bf Acknowledgments:} V.M.~thanks Gerardo Morsella for stimulating discussions and {Karl-Hermann Neeb for useful comments}.
V.M. acknowledges the support of  the University of Rome Tor Vergata funding \emph{OANGQS} CUP E83C25000580005,  the MUR Excellence Department Project \emph{MatMod@TOV} awarded to the Department of Mathematics, University of Rome Tor Vergata, CUP E83C23000330006 and INdAM-GNAMPA.

\medskip

\nin {\bf Notation} 

\begin{itemize}
\item The neutral element of a group $G$ is denoted $e$, and
  $G_e$ is the identity component.
\item The Lie algebra of a Lie group $G$ is denoted $\L(G)$ or $\g$.   
\item For an involutive automorphism $\sigma$ of $G$, we write
  $G^\sigma =  \{ g\in G \:  \sigma(g) = g \}$ for the subgroup of fixed points
  and $G_\sigma := G \rtimes \{\id_G, \sigma\}$ for the corresponding
 group extension. 
\item $\AU(\cH)$ is the group  of unitary or antiunitary operators
    on a complex Hilbert space. 
\item $\cC_0^\infty(M)$ denotes the set of compactly supported functions in $M$.
\item An (anti-)unitary representation of $G_\sigma$ is a homomorphism
    $U \: G_\sigma \to \AU(\cH)$ with $U(G) \subeq \U(\cH)$ for which 
    $J := U(\sigma)$ is antiunitary, i.e., a conjugation.
\item Unitary or (anti-)unitary representations on the complex Hilbert space
  $\cH$ are denoted as pairs $(U,\cH)$. For a unitary representation $(U,\cH)$ of a Lie group $G$ we write:
$\partial U(x) = \derat0 U(\exp tx)$ for the infinitesimal
    generator of the unitary one--parameter group $(U(\exp tx))_{t \in\R}$
    in the sense of Stone's Theorem. In this setting $C_U := \{ x \in \g \: -i \partial U(x) \geq 0\}$.
  \item An element $h$ of a Lie algebra $\g$ is called
    \begin{itemize}
    \item     {\it hyperbolic} if $\ad h$ is diagonalizable over $\R$
    \item {\it elliptic} or {\it compact} if $\ad h$ is semisimple
      with purely imaginary spectrum.
    \end{itemize}
    \item Let $\cA\subset \cB(\cH)$, then the commutant of $\cA$ is $\cA'=\{b\in\cB(\cH):ab=ba, \forall a\in\cA\}$.	\\ A  (complex) $*$-subalgebra  $\cA\subset\cB(\cH)$ is a von Neumann algebra if and only if $\cA=\cA''$. A vector $\Omega\in\cH$ is cyclic for $\cA$ if $\overline{\cA\Omega}=\cH$ and it is separating if $\overline{\cA'\Omega}=\cH.$

    \item $\dd U \colon \g \to \End(\cH^\infty)$ for the representation of
    the Lie algebra $\g$ on the space $\cH^\infty$ of smooth vectors. Then
    $\partial U(x) = \oline{\dd U(x)}$ (operator closure) for $x \in \g$. 
\end{itemize}

\section{Preliminaries}
\mlabel{sec:2}
In this section we introduce the notations, definitions,  theorems and the general techniques that will be employed throughout the paper. We will start by examining the geometry of fundamental examples in Section \ref{sect:wre}. From them it will emerge the key role of Euler elements and the structure of Euler wedges as explained in Section \ref{sect:ee}. Wedge domains in causal homogeneous spaces are briefly recalled in Section \ref{sect:chs}. 
In Sect.~\ref{sect:geoss} we present the fundamental free field construction (Sect.~\ref{sect:ff}) and we introduce the standard subspace formalism (Sect.~\ref{sect:stsub}) to focus on the axiomatic framework for nets of standard subspaces on abstract wedges and on the generalized BGL-construction provided in \cite{MN21} (Sect.~\ref{sect:BGL}).

\subsection{Wedge regions in some spacetimes}\label{sect:wre}

In this section we introduce the wedge regions in some of the spacetimes that are relevant for physics and state the correspondence between these regions, specific subgroups of the symmetry group and Lie algebra elements.

Consider  Minkowski spacetime $M=\RR^{1,d}$  of space dimension $d$. Here the metric considered is
\[ ds^2=dx_0^2-dx_1^2-\ldots-dx_d^2. \]
Let $x^2$ be the square of the Minkowski length, one defines \textit{timelike vectors} by $x^2>0$, \textit{lightlike vectors or null vectors} by $x^2=0$, and \textit{spacelike vectors} by $x^2<0$. We shall denote $C=\{x\in M:x^2\geq0, x_0\geq0\}$ the pointed convex closed cone describing the timelike and lightlike future of the origin.
Given a region $\cO$, its causal complement is $$\displaystyle{\cO'=\{y\in\RR^{1,d}: (x-y)^2<0, \forall x\in\cO\}^\circ}.$$ $\cO' $  represents the set of points that are causally independent of those in $\cO$. A causally complete region satisfies the property $\cO=\cO''$.
The connected component of the identity of the homogenous symmetry group of $M$ is the proper orthochronous Lorentz group $\cL_+^\up=\SO(1,d)_e$. We denote by $\cL_+$ the proper Lorentz group, namely the group generated by $\cL_+^\up$ and the space and time reflection  \[ j_{W_R}(x) =(-x_0,-x_1,x_2,\ldots,x_d) \]  acting on the first two coordinates, hence $\cL_+=\cL_+^\up\cup\cL_+^\downarrow$ where $\cL_+^\downarrow=j_1\cL_+^\up$. 
The associated inhomogeneous symmetry groups  denoted by $\cP_+^\up=\RR^{1,d}\rtimes \cL_+^\up$ and $\cP_+=\RR^{1,d}\rtimes \cL_+$ are respectively the proper orthochronous Poincar\'e group and the proper Poincar\'e group. In particular $\cP_+=\cP_+^\up\cup\cP_+^\downarrow=\cP_+^\up\cup j_1\cP_+^\up$.

In this framework there are some distinguished regions (connected open subsets) of $M$ that have a deep connection with the symmetry group.
A \textit{wedge region }$$W=gW_R\subset\RR^{1,d}\quad g\in \cP_+^\up$$ is a Poincar\'e transform of the  Rindler wedge  $W_R$ (here also called the right wedge) introduced in~\eqref{eq:Rinwed}.

Now, consider  the one--parameter group of boosts $\Lambda_{W_R}(t)$ contained in $\cP_+^\uparrow$, as in \eqref{eq:boost}.
All the other one-parameter group of boosts are determined by the Poincar\'e adjoint action $\Lambda_W=g\Lambda_{W_R}g^{-1}$, with $ g\in\cP_+^\up$. 
There is a 1-1 correspondence between wedge regions and one-parameter group of boosts and their generators as
$$h_W\stackrel{1:1}\longleftrightarrow  \Lambda_W=g\Lambda_{W_R}g^{-1}\stackrel{1:1}\longleftrightarrow W=gW,\qquad  g \in\cP_+^\up, $$
where $h_W\in\Lie(\cP_+^\up)$ generates $\Lambda_W(t)=\exp(th_W)$. 

Let $W=gW_R$ and $W'=\{x\in M: (x-y)^2<0,y\in W\}^\circ$, the transformation $j_W:=gj_{W_R}g^{-1}\in \cP_+^\down$  maps $W$ onto $W'$. Clearly $j_W=j_{W'}$ and it is  important to remark that $j_W=\exp(i\pi h_W)$.
We have observed that the \textit{wedge regions can be characterized by Lie algebra elements} and that there exists a natural transformation mapping $W$ to $W'$, which coincides with the analytic extension of the one-parameter group $\Lambda_W(t)$. We further remark  that $\Lambda_{W'}(t)=\Lambda_W(-t)$, hence $h_{W'}=-h_W$. 
When conformal symmetries are considered - specifically, when the symmetry group consists of local diffeomorphisms (defined outside meager sets) that preserve the metric tensor up to a positive function - all conformal transformations of a wedge in Minkowski space, such as wedges, double cones, future cones, and past cones, are referred to as conformal wedge, see \cite{BGL02,BGL93}.  Here, a future and a past cone are intended to be respectively  regions of the form $(C+a)^\circ$ and $(-C+b)^\circ$ for some $a,b\in M$ and a double cone is a region of the form  $(C+a)^\circ\cap(-C+b)^\circ$ where $b-a\in C^\circ$. The identification of conformal wedges with conformal Lie generators remains valid.

Consider the $d-$dimensional de Sitter spacetime as the subset $$\dS^d=\{x\in\RR^{1,d}:x^2=-1\}\subset \RR^{1,d},$$ endowed with the metric obtained by restriction of  the Minkowski metric. The symmetry group of isometries is the Lorentz group $\cL_+^\up$ and its $\mathbb Z_2$-graded extension $\cL_+$.  Wedges are defined as $W^{\rm dS}:=W\cap \dS^d$, where $W=gW_R$ is a Minkowksi wedge, with $g\in \cL_+^\up$. The one--parameter groups of boosts are $\Lambda_{W^{\dS}}:=g\Lambda_{W_R^{\dS}}g^{-1}$  with $g\in \cL_+^\up$, where $\Lambda_{W_R^{\dS}}=\Lambda_{W_R}\in\cL_+^\up$ and there is a 1-1 correspondence 
$$ h_{W^{\dS}}\stackrel{1:1}\longleftrightarrow\Lambda_{W^{\dS}}(t)=\exp(th_{W^{\dS}})\stackrel{1:1}\longleftrightarrow W^{\dS}$$ where $h_{W^{\dS}}\in \Lie (\mathcal{L}_+^\up)$ and $\exp(h_{W^{\dS}}t)=\Lambda_{W^{\dS}}(t)$.
As in the Poincar\'e case one can define the causal complement of a wedge region ${W^{\dS}}'$  and the reflection $j_{W^{\dS}}$ which satisfies  $j_{W^{\dS}}=\exp(i\pi h_{W^{\dS}})$ and $j_{W^{\dS}}{W^{\dS}}={W^{\dS}}'$.

As a last example we  consider the chiral circle $S^1$.  The symmetry group is the M\"obius group $\Mob$, the group of orientation preserving fractional linear transformations fixing $S^1\subset\mathbb C$. One can pass from  this picture to the one-point compactification of the real line $\bar\RR=\RR\cup\{\infty\}$ through the Cayley transform $\bar\RR\ni x\mapsto -\frac{x-i}{x+i}\in S^1$. Here the symmetry group becomes $\Mob\simeq \PSL_2(\RR)$.
On $S^1$ the notion of wedge regions is attached to connected, open, non-empty, non dense intervals $I$ of the circle.  In the real line picture $\bar \RR$, they correspond to bounded intervals, half-lines or the complement of a bounded interval of $\RR$ in $\bar \RR$.
There is a 1-1 correspondence  between  one--parameter groups of ``boost"-dilations in  $\Mob$ and ``wedge"-intervals in $S^1$ given by 
$$h_I\in\Lie(\Mob)\stackrel{1:1}\longleftrightarrow\delta_I(t)=\exp(h_I t)\in\Mob\stackrel{1:1}\longleftrightarrow I\subset S^1$$ 
where $\delta_{I=gI_0}(t)=g\delta_{I_0}(t)g^{-1}$ with $g\in \Mob$ and $\delta_{I_0}(t)$ in the real line picture given by $\delta_{I_0}(t)x=e^t x$ for $I_0=(0,+\infty)$.
In this picture the wedge reflection $j_I=\exp(i\pi h_I)$ is  sending $I$ to $I':=\Int(S^1\backslash I)$ reversing the point ordering.
Also in this case $\delta_{I'}(t)=\delta_I(-t)$, hence $h_{I'}=-h_I$, and $j_{I'}=j_I$. 
This example becomes relevant when one considers the two dimensional Minkowksi space with chiral coordinates $(u,v)=({x_0+x_1},{x_0-x_1})$ with  $\Mob\times \Mob$ as the symmetry group, cf.~\cite{BGL93,MT19}.

We  will now observe that these examples of correspondences between wedges and one-parameter groups of boost generators are special cases of a more general framework.

\subsection{Euler elements and abstract wedges}\label{sect:ee}
 
In all the previous examples any element $h_W$ identifying a wedge region $W$ in the related space-time $M$ has the property that the adjoint action of the Lie element on the Lie algebra is diagonalizable with eigenvalues $\{-1,0,1\}$, namely $\Spec(\ad h_W) \subeq \{-1,0,1\}$. This is particularly evident for the 2+1 dimensional de Sitter spacetime, where $\cL_+=\PSL_2(\RR)$, and $\Lie(\cL_+^\up)\simeq\fsl_2(\RR)$
 and 
 ${\displaystyle h_{W^{\dS}_R}=\frac12\left(\begin{array}{cc}1&0\\0&-1\end{array}\right)}$ provides the standard 3-grading of $\fsl_2(\RR)$:

$$\mathfrak{l}=\mathfrak{l}_{-1}\oplus\mathfrak{l}_0\oplus \mathfrak{l}_{-1}=\RR\left(\begin{array}{cc}0&1\\0&0\end{array}\right)\oplus\RR\left(\begin{array}{cc}1&0\\0&-1\end{array}\right)\oplus\RR \left(\begin{array}{cc}0&0\\1&0\end{array}\right).$$ 

The same holds for the M\"obius group  $\Mob\simeq \PSL_2(\RR)$ in the real line picture.
On the other hand, there is a difference in the wedge geometry:  de Sitter spacetime does not admit proper wedge inclusions, unlike the conformal circle, where wedge inclusions correspond to interval inclusions.

 We now introduce the natural general framework in which all these examples fit, clarifying the previously mentioned differences regarding the inclusion properties of wedge regions. The primary reference for this is \cite{MN21}. 

An element $h$ in a finite dimensional 
real Lie algebra $\g$ is called an
{\it Euler element} if $\ad h$ is non-zero and diagonalizable with 
$\Spec(\ad h) \subeq \{-1,0,1\}$. In particular the eigenspace 
decomposition with respect to $\ad h$ defines a $3$-grading 
of~$\g$: 
\[ \g = \g_1(h) \oplus \g_0(h) \oplus \g_{-1}(h), \quad \mbox{ where } \quad 
\g_\nu(h) = \ker(\ad h - \nu \id_\g)\] 
Then {$\tau_h(y_j) = \exp (i\pi \ad h)(y_j)= (-1)^j y_j$} for $y_j \in \g_j(h)$ 
defines an involutive automorphism of $\g$  called {an {\it Euler involution on} $\fg$}. 

We denote $\cE(\g)$  the set of Euler elements in~$\g$. 
The orbit of an Euler element  $h$ under the group 
$\Inn(\g) = \la e^{\ad \g} \ra$ of 
{\it inner automorphisms} 
is denoted with $\cO_h = \Inn(\g)h \subeq \g$.
We say that $h$ is {\it symmetric} if $-h \in \cO_h$. 
In a simple Lie algebra $\fg$, the symmetric Euler elements are always contained in a subalgebra $\fh\simeq\fsl_2(\RR)$  and non-symmetric Euler elements are contained in a subalgebra $\fh\simeq\fgl_2(\RR)$ (see \cite[Lem.~2.15]{MN22} for a more general case).
Let $G$ be a connected Lie group with Lie algebra $\fg$ and $\cE(\fg)\neq\emptyset$, {assume that $\tau$ is the involution on $G$ integrating an Euler involution on $\fg$.} Then $G_{\tau}=G\rtimes_{\tau} \mathbb Z_2=G\rtimes\{e,\tau\}$ is well-defined and one has a continuous homomorphism
$\epsilon \:  G_{\tau} \to \{\pm 1\}$ where
\[ G = \epsilon^{-1}(1) \quad \mbox{ and } \quad  G\,\tau = \epsilon^{-1}(-1),\] 
so that $G \trile G_{\tau}$ is a normal subgroup of index $2$ and $G\tau = G_{\tau} \setminus G$. 
Note that if $G= \Inn(\fg)$ or $G$ is simply connected, then $\tau_h$ always integrate to an involution on $G$ when $h\in\cE(\fg)$. 

The set 
{  \[ \cG_E {:= \cG_E({G_\tau})} :={ \{ (h,\sigma)\in\g \times G\tau \:
      \, h\in\cE(\fg), \Ad\sigma=\tau_h\}}\] 
 is called the {\it abstract Euler wedge space} of ${G_\tau}$, an element $(h,\sigma) \in \cG_E$ is an {\it Euler couple} or {\it Euler wedge} and $\sigma$ is called \textit{Euler involution}. For sake of simplicity we here assume here that that  $G$ is center free. In this case every Euler element identifies an Euler couple $(h,\tau_h)\in\cG_E$ uniquely. }
One can also define more general abstract wedges without referring to Euler elements, see \cite[Sect.~2]{MN21}.

Given a fixed couple $W_0 = (h,\tau_h) \in \cG_E$, the orbits 
\[ \cW_+ (W_0):= G.W_0 \subeq \cG_E  \quad \mbox{ and } \quad 
\cW(W_0) := G_{\tau}.W_0 \subeq \cG_E  \] 
are called the {\it positive} and the {\it full
 {abstract} wedge space containing $W_0$}. 

\textbf{Covariance of Euler wedges.} A $G_{\tau}$-covariant action on $\cG_E $ can be defined as follows.
Consider the 
{\it twisted adjoint action} of $G_{\tau}$ which changes the sign on odd group elements: 
\begin{equation}
  \label{eq:adeps}
  \Ad^\epsilon \: G_{\tau} \to \Aut(\g), \qquad 
\Ad^\epsilon(g) := \epsilon(g) \Ad(g).
\end{equation}
Then $G_{\tau}$ acts on $\cG_E$  by 
\begin{equation}
  \label{eq:cG-act}
 g.(h,\tau_h) := (\Ad^\epsilon(g)h, g\tau_h g^{-1}).
\end{equation}

\textbf{Order structure on Euler wedges.} Given an $\Ad^\epsilon(G_{\tau})$-invariant pointed closed convex cone $C_\fg \subeq \fg$, consider the subsets $$C_\pm := \pm C_\g \cap \g_{\pm 1}(h),$$
and the compression semigroup,
\[ S_W := \exp(C_+) G_{W} \exp(C_-) 
  = G_{W} \exp\big(C_+ + C_-\big)\] where $G_W := \{ g \in G \:  g.W = W\}$, see \cite[Thm.~2.16]{Ne22}.
Then the partial order on the orbit 
$G.W \subeq \cG_E$ is defined by 
\begin{equation}
  \label{eq:cG-ord}
g_1.W \leq g_2.W  \quad :\Longleftrightarrow \quad 
g_2^{-1}g_1 \in S_W.
\end{equation}
In particular, $g.W \leq W$ is equivalent to $g \in S_W$. 
As an example, consider $G=\PSL_2(\RR)$ and its $\mathbb Z_2$-graded extension $G_{\tau_h}$.  If 
$C$ is chosen to be trivial, we arrive at a situation in which no non-trivial wedge inclusions exist, corresponding to the de Sitter wedge geometry. If we choose 
\begin{equation}\label{eq:Csl}
C: 
= \Big\{  X= \pmat{a & b \\ c & -a} \: 
b \geq 0, c \leq 0, a^2 \leq -bc\Big\}\subset\fsl_2(\RR) 
\end{equation}
we are in the situation where there are proper wedge inclusions, hence the case of the circle $S^1$ wedge geometry.
Namely the choice of $C$ influences the geometry of the model.

\textbf{Duality and Locality properties of  Euler wedges.} 
For $W = (h,\tau_h) \in \cG_E$, we define the {\it dual wedge} by   
$W' := (-h,\tau_h) {= \tau_h.W}$. 
This definition aligns with the geometric relation of a wedge with its causal complement in spacetime manifolds.
Note that $(W')' = W$ and $(gW)' = gW'$ for $g \in G$ 
by~\eqref{eq:cG-act}.  
Since the wedge inclusions and the dual wedge are well defined notions on couples of Euler wedges we can state the \textit{locality property of Euler couples} $W_1$ and $W_2$ that satisfy $W_1\subset W_2'$. 

 Note that it is not always the case that $W_0'\in \cW_+(W_0)$. As an example, consider the translation-dilation group $G$ acting on the real line. A wedge is a right half-line $I=\RR^++a$ and is associated 1-1 to the positive dilations $\delta_I(t)x=\exp(h_It)=e^{t}(x-a)+a$ where $x\in\RR$. The complement interval  $I'$ is a negative half-line $I'=\RR^-+a$, it has associated the one-parameter group $\delta_{I'}(t)=\delta_I(-t)$. There exists no $g\in G$ such that $gI=I'$ and one needs to include the $\tau_I$ to reflect $I$ to $I'$, namely $I'\in\cW(I)$ but $I'\not\in\cW_+(I)$. This correspond to the fact that $h$ is not symmetric in $\fg$, see  \cite{MN21}.
We further remark that if  $Z(G)$ is  non trivial, many abstract wedges correspond to an Euler element, and one can define a notion of ``twisted  complements" for a more general notion of complement wedges. In this case it is not always the case that $W'_0\in\cW_+(W_0)$ even if $h$ is symmetric. This case is discussed in \cite[Sect.~2.4.2]{MN21}.

We have now established a suitable framework for defining an abstract index set for Quantum Field Theory. Notably, this structure extends beyond well-known examples within the physics context. One can classify the real Lie algebras  $\fg$  supporting a three-grading given by Euler elements. In particular, the Lie group $\Inn(\fg)$, hence its extension $\Inn(\fg)_{\tau_h}$, supports Euler wedges. This generalizes the fundamental examples outlined in the previous section.
By Theorem 3.10 in \cite{MN21} if a simple real Lie algebra $\fg$ supports Euler elements then its restricted root system has to be one of the following

\begin{multicols}{2}
\begin{itemize}
\item $A_n: n\,\Inn(\fg)\text{-orbits  of Euler elements} $
\item $ B_n: \text{one} \Inn(\fg)\text{-orbit  of Euler elements} $
\item $ C_n: \text{one}\Inn(\fg)\text{-orbit of Euler elements} $
\item $ D_n:  \text{three} \Inn(\fg)\text{-orbits  of Euler elements}  $
\item $E_6: \text{two}\Inn(\fg)\text{-orbits  of Euler elements} $
\item $E_7:  \text{one}\Inn(\fg)\text{-orbit  of Euler elements} $
\end{itemize}
\end{multicols}
The theorem also classify the number of  inequivalent $\Inn (\fg)$-orbits.
Here is the list of three graded real Lie algebras:
\medskip

\begin{tabular}{||l|l|l|l|l||}\hline
& $\g$  & $\Sigma(\g,\fa)$  & $h$ & $\g_1(h)$  \\ 
\hline\hline 
1 & $\fsl_n(\R)$ & $A_{n-1}$ & $h_j, 1 \leq j \leq n-1$ & $M_{j,n-j}(\R)$  \\ 
2 & $\fsl_n(\H)$ & $A_{n-1}$ & $h_j, 1 \leq j \leq n-1$ & $M_{j,n-j}(\H)$  \\ 
3 & $\su_{n,n}(\C)$ & $C_{n}$ & $h_n$ & $\Herm_n(\C)$  \\ 
4 & $\sp_{2n}(\R)$ & $C_{n}$ & $h_n$ & $\Sym_n(\R)$   \\ 
5 & $\fu_{n,n}(\H)$ & $C_{n}$ & $h_n$ & $\Aherm_n(\H)$  \\ 
6  & $\so_{p,q}(\R)$ & $ B_p\ (p<q),\ D_p\ (p = q)$ & $h_1$ & $\R^{p+q-2}$   \\ 
7  & $\so^*(4n)$ & $C_n$ & $h_n$ & $\Herm_n(\H)$   \\ 
8  & $\so_{n,n}(\R)$ & $C_n$ & $h_n$ & $\Alt_n(\R)$   \\ 
9 & $\fe_6(\R)$ & $E_6$ & $h_1=h_6' $ & $M_{1,2}(\bO_{\rm split})$   \\ 
10 & $\fe_{6(-26)}$ & $A_2$ & $h_1$ & $M_{1,2}(\bO)$   \\ 
11 & $\fe_7(\R)$ & $E_7$ & $h_7 $ & $\Herm_3(\bO_{\rm split})$   \\ 
12 & $\fe_{7(-25)}$ & $C_3$ & $h_3$ & $\Herm_3(\bO)$   \\ 
13 & $\fsl_n(\C)$ & $A_{n-1}$ & $h_j, 1 \leq j \leq n-1$ & $M_{j,n-j}(\C)$  \\ 
14 & $\sp_{2n}(\C)$ & $C_{n}$ & $h_n$ & $\Sym_n(\C)$   \\ 
15a & $\so_{2n+1}(\C)$ & $ B_{n}$ & $h_1$ & $\C^n$   \\ 
15b & $\so_{2n}(\C)$ & $ D_{n}$ & $h_1$ & $\C^n$   \\ 
16 & $\so_{2n}(\C)$ & $ D_{n}$ & $h_{n-1}, h_n$ & $\Alt_n(\C)$   \\ 
17 & $\fe_6(\C)$ & $E_6$ & $h_1 = h_6' $ & $M_{1,2}(\bO)_\C$   \\ 
18 & $\fe_7(\C)$ & $E_7$ & $h_7 $ & $\Herm_3(\bO)_\C$   \\ 
\hline
\end{tabular} \\[2mm] {\rm Table 1: Simple $3$-graded Lie algebras. We follow the conventions of the tables in [Bo90a] for the classification of irreducible root systems, the enumeration of the simple roots and the associated Euler elements, cf.~\cite{MN21}}\\

One can further classify the orbits of symmetric Euler elements, leading to a reduced classification as follows:
\begin{itemize}
\item $A_{2n-1}: \text{one} \Inn(\fg)\text{-orbit of symmetric Euler elements}$
\item $B_n: \text{one} \Inn(\fg)\text{-orbit of symmetric Euler elements} $
\item $C_n: \text{one} \Inn(\fg)\text{-orbit of symmetric Euler elements} $
\item $D_{2n+1}: \text{one} \Inn(\fg)\text{-orbit of symmetric Euler elements} $
\item $D_{2n}: \text{three}  \Inn(\fg)\text{-orbits of symmetric Euler elements} $
\item $E_7:\text{one} \Inn(\fg)\text{-orbit of symmetric Euler elements}$
\end{itemize}

\noindent{Here is the list of simple hermitian Lie algebras $\g^\circ$ 
supporting Euler elements. Theorem \cite[Prop.~3.11]{MN21} ensures that they have to be of tube type, namely their restricted root system is of type $C_r$.}

\medskip

\begin{tabular}{||l|l|l|l|l||}\hline \label{tab:2}
{} $\g^\circ$ \mbox{(hermitian)}  & $\Sigma(\g^\circ, \fa^\circ)$ & $\g = (\g^\circ)_\C$ & $\Sigma(\g,\fa)$ & {\rm symm.\ Euler element}\  $h$ \\ 
\hline\hline 
$\su_{n,n}(\C)$ & $C_n$ & $\fsl_{2n}(\C)$ & $A_{2n-1}$ & $h_n$ \\ 
 $\so_{2,2n-1}(\R), n > 1$ & $C_2$ & $\so_{2n+1}(\C)$ & $B_n$ & $h_1$ \\ 
$\sp_{2n}(\R)$ & $C_n$ & $\fsp_{2n}(\C)$ & $C_n$ & $h_n$ \\ 
 $\so_{2,2n-2}(\R), n > 2$ & $C_2$ & $\so_{2n}(\C)$ & $D_{n}$ & $h_1$ \\ 
 $\so^*(4n)$ & $C_n$ & $\so_{4n}(\C)$  & $D_{2n}$ & $h_{2n-1},  h_{2n}$ \\ 
$\fe_{7(-25)}$ & $C_3$ & $\fe_7$ & $E_7$ & $h_7$ \\ 
\hline
\end{tabular} \\[2mm] {\rm Table 2: Simple hermitian Lie algebras $\g^\circ$ 
of tube type}\\

 \subsection{Wedge domains in causal homogeneous spaces} 
 \label{sect:chs}

{Let $G$ be a connected Lie group, $H$ a closed subgroup. The associated  homogeneous  space is given by the quotient $M=G/H$. }It is said to be \textit{causal} if there exists a $G$-invariant field of pointed generating closed convex cones $\{C_m\}_{m\in M}$ with $C_m\subset T_m(M)$ the tangent space at the point $m\in M$.  The Minkowski space, as well as the de Sitter space, is  a causal homogeneous space, because of the causal structure given by the Minkowksi metric. 

{Given a causal manifold $M$ and an Euler element $h\in\fg$,  the {\it modular vector field} is

\begin{equation}
  \label{eq:xhdef}
 X_h^M(m) 
 :=  \frac{d}{dt}\Big|_{t = 0} \exp(th).m .
\end{equation}
Then the connected component
\begin{equation}
  \label{eq:W-def}
  W := W_M^+(h)_{eH}
\end{equation}
of the base point
$eH \in M$ in the {\it positivity region}
\begin{equation}\label{def:WM}
 W_M^+(h) := \{ m \in M \: X^M_h(m) \in C_m^\circ \} 
 \end{equation}
 is the natural candidate for a wedge domain if $eH\in W_M^+(h)$ (or its boundary). Motivated by the Bisognano--Wichmann property (BW) and its consequences in AQFT,  the modular flow should  be indeed timelike future-oriented in the wedge region. 
It  corresponds to the inner  time evolution of Rindler wedge (see \cite{CR94} and also
\cite{BB99, BMS01}, \cite[\S 3]{CLRR22}). As a concrete example, in Minkowski and de Sitter spacetimes the wedge $W_\Lambda$ is determined as the set of points $x\in M$ such that  the vector field at $x$ associated to the one-parameter group $\Lambda$  of boosts is timelike and future pointing.}

{In this context we shall consider causal symmetric spaces that are homogeneous space of the form $M = G/H$, where $G$ is a connected Lie group, and there exists an involution $\tau$ on $G$ for which $H \subset G^\tau$  is a union of connected components. We can assume $G$ to be simple and center free. 
Let $C\subset \fg^{-\tau}\simeq T_{eH}(M)$ be a pointed generating closed $H$-invariant convex cone producing the causal structure on $M$ by $G$-action. We can distinguish symmetric spaces that are \textit{compactly causal}, namely the cone $C^\circ$ contains only elliptic elements and \textit{non-compactly causal symmetric spaces}, namely the cone $C^\circ$ contains only hyperbolic elements} (see \cite{HO96}).
The Anti-de Sitter spacetime belongs to the first family and  de Sitter spacetime belongs to second one.
The wedge regions, as introduced above, have been studied for compactly and
non-compactly causal symmetric spaces in \cite{NO23} and \cite{NO22,MNO23a,MNO23b}, respectively.

Assume that the center of $G$ is trivial, then there is a 1-1 correspondence between Euler elements and abstract Euler wedges (Euler couples).  In particular the stabilizers $G^{(h,\tau_h)}$ and $G^h$ coincide
    and one can identify a wedge orbit $\cW_+(h,\tau_h) \subeq \cG_E$ with the adjoint orbit $\cO_h = \Ad(G)h$. We thus obtain a natural map
  from $\cW_+(h,\tau_h) \cong \cO_h$ to regions in $M$
 by $g.(h,\tau_h) \mapsto g.W_M^+(h)$. If, in addition, 
  $G^h$ preserves the connected component $W \subeq W_M^+(h)$, as it happens if $W_M^+(h)$ is connected, one can define a map from the abstract wedge space $\cW_+(h,\tau_h) \ni W=(k,\tau_k)$ to the
  geometric wedge regions $W_M^+(h)$ on~$M$.  

{Given a real connected simple Lie group $G$ such that the Lie algebra $\fg$ contains an Euler element $h\in\cE(\fg)$, one can construct a non-compactly causal symmetric space as $M=G/H$ where $H$ is an open subgroup of $G^{\tau_{nc}}$, the centralizer of the involution $\tau_{nc}$ that integrates $\tau=\tau_h \theta$  on $G$ where $\theta$ is a  Cartan involution of $\fg$ such that $\theta (h)=-h$. Note that within this picture $h\in T_\textbf{e}(M)\simeq\fg^{-\tau}$ where $\textbf{e}=eH$ is the basepoint of $M$. Furthermore one can define a $G$-invariant field of pointed generating closed convex cones $C_\textbf{p} \subset T_\textbf{p}(M)$ such that $h\in C_\textbf{e}^\circ$, hence $\textbf{e}\in W_M^+(h)$, and $C_\textbf{e}$ is $H$-invariant, see \cite{MNO23a}. In this way one can associate a wedge domain $W(h)\subset M$ to any $h\in\cE(\fg)$. A  classification result of the local structure of such non compactly-causal homogeneous spaces and connectedness of $W^+_M(h)$ is described in \cite[Thm.~4.21]{MNO23a} and \cite[Thm.~7.1]{MNO23b}. }

As an example, consider the de Sitter space $\dS$ of dimension $d$ and its symmetry group, the Lorentz group, $G=\SO(1,d)_e$. Let $h:=h_{W^{\dS}_R}$ be the generator of the one-parameter group of boosts in \eqref{eq:boost}. Then following \cite{MNO23a},  one can consider the involution $\tau=\tau_{h}\theta$ on the Lie algebra $\fg=\so(1,d)$ where $\theta$ is the Cartan involution of the Lie algebra $\fg$ promoted to $G$. In particular $$\theta=\Ad\left(\diag(-1,\textbf{1}_{d})\right),\quad\tau_h=\Ad\left(\diag(-1,-1,\textbf{1}_{d-1})\right),\quad {\tau=\Ad\left(\diag(1,-1,\textbf{1}_{d-1})\right)}.$$
Let $\textbf{e}_1=(0,1,0,\ldots,0)\in \dS\subset \RR^{1,d}$ and $G_{\textbf{e}_1}$ be its stabilizer group, then the centralizer of $\tau_{nc}$ in $G$ is described explicitly by 
{$$G^{\tau_{nc}}\simeq G_{\textbf{e}_1}\rtimes \{1,r_{1,2}(\pi)\}\simeq G_{\textbf{e}_1}\rtimes \Z_2$$ }where $r_{1,2}(\pi)$ is the $\pi$-rotation in the first two {spatial} coordinates. In particular $G^{\tau_{nc}}$ has two connected components and $G_{\textbf{e}_1} $ is the identity component fixing $\textbf{e}_1$. Then   de Sitter space $\dS$ is isomorphic to $G/G_{\textbf{e}_1}$ which is a causal symmetric space but  the quotient $G/G^{\tau_{nc}}$ is the projective de Sitter space  which is not a causal space. With the identification $\dS\simeq G/G_{\textbf{e}_1}$ one has that $W_M^+(h)=W^{\dS}_R$. 
In general, let $\fg$ be a non compactly causal simple Lie algebra, $G=\Inn(\fg)$ and $h\in\cE(\fg)$, then $M=G/G^{\tau_{nc}}$ is a causal  if and only if $h$ is not symmetric. Details are in \cite[Thm.~4.19]{MNO23a}).  

\subsection{ The geometry of nets of real subspaces}
\label{sect:geoss}

In this section we recall some fundamental properties of the
{geometry of} the standard subspaces and the one-particle nets.
{We refer to \cite{LRT78, BGL02, MN21, NO17, Lo08} for more details}.
Particles in Quantum Field Theory are defined  through their symmetries rather than pointwise properties and their meaning is clear when the symmetry group is large enough. For instance, one-particle states of a free Quantum Field Theory on Minkowski spacetime   are unit vectors of the Hilbert space supporting a unitary irreducible positive energy representation of the Poincar\'e group.

The free field constitutes a fundamental model in AQFT. We will see how the language of the standard subspaces captures  many of the features of the free models and how this   allows a generalization of the AQFT framework on abstract sets of wedges and on causal symmetric spaces.

\subsubsection{Scalar free fields in Algebraic Quantum Field Theory} \label{sect:ff}

The free scalar mass $m\geq0$ Quantum Field Theory on Minkowksi spacetime can be constructed as follows. Let $U$ be the unitary, scalar, mass $m$ positive energy representation of the proper  Poincar\'e group $\cP_+^\up=\RR^{1,3}\rtimes \SO(1,3)_e$.  The Hilbert space
is $\cH=L^2(\Omega_m,d\Omega_m)$, where $$\Omega_m=~\{p\in~M:~p^2=m^2, p_0>0\}$$ and $d\Omega_m$ is the Lorentz invariant measure on $\Omega_m$. 
Let $f\in~\cS(M)$, define $\hat f(p)=\frac1{(2\pi)^{2}}\int_M e^{ipx}f(x)dx$.  Then the functions $$\hat \cS(M)=\{\hat f|_{\Omega_m}:f\in\cS(M)\}$$ defines a complex linear dense subset of $\cH$. 
On this Hilbert space the unitary scalar representation acts by
$$(U(g)f)(x)=f(g^{-1}x),\qquad g\in \cP_+^\up.$$ Positivity of the energy means that the joint spectrum of the translation generators is contained in the forward lightcone $C$. This representation extends { by the same formula }to a (anti-)unitary representation of the proper Poincar\'e group $\cP_+$, \cite[Thm.~9.10]{Var85}.
To any open bounded $\cO\subset M$ we can associate  the real subspace
$$\sH(\cO)=\overline{\{\hat f|_{\Omega_m}:f\in\cC_0^\infty(M), \supp f\,\subset \cO\}}.$$ This is the set of  one-particle states localized in $\cO$. Note that the following properties  hold
\begin{itemize}
\item \textbf{Isotony:} if $\cO_1\subset\cO_2$ then $\sH(\cO_1)\subset \sH(\cO_2)$,
\item \textbf{Locality:} if $\cO_1\subset\cO_2'$, then $\sH(\cO_1)\subset \sH(\cO_2)'$, where $\sH(\cO_2)'$ is the symplectic complement of $\sH(\cO_2)$ defined by $\sH(\cO_2)' :=\{\xi\in\cH:\  \Im\langle\xi,\eta \rangle=0, \forall \eta \in \sH(\cO_2)\} $,
\item \textbf{Covariance:} if $g\in\cP_+$, then $U(g) \sH(\cO)= \sH(gO)$,
\item \textbf{Reeh--Schlieder:} $\overline{\sH(\cO)+i\sH(\cO)}=\cH$ for all $\cO\neq\emptyset$.
\end{itemize}
Wedge subspaces can be defined as 
\begin{equation}\label{eq:ffn}\sH(W)=\overline{\sum_{\cO\subset W}\sH(\cO)}
\end{equation}
The free field is the second quantization net of von Neumann algebras on the Fock Hilbert space. It is constructed as follows. Given the one-particle Hilbert space $\cH$, we can define the Fock space $$\cF(\cH)=\bigoplus_{n=1}^\infty\cH^{\otimes_s^n}\oplus \C\cdot\Omega,$$
where $\cH^{\otimes_s^n}$ is the symmetrized $n$-fold tensor product of the one-particle Hilbert space $\cH$. 
For every vector $\xi\in\cH$ the associated Weyl unitary $\mathrm{w}(\xi)$ on the Fock space is defined by the commutation relations
\begin{equation}\label{eq:weyl}\mathrm{w}(\xi)\mathrm{w}(\eta)=e^{-\frac12\Im(\xi,\eta)}\mathrm{w}(\xi+\eta),\qquad \xi,\eta\in\cH\end{equation}
and the expectation value on the vacuum state
\begin{equation}\label{eq:vs}(\Omega,\mathrm w(\xi)\Omega)=e^{-\frac14\|\xi\|^2}.\end{equation}  Then the map $\cH\ni h\mapsto \mathrm{w}(h)\in\cF(\cH)$ is strongly continuous. Further details in \cite{BR87}.
Given a real subspace $H\subset \cH$, we can associate a von Neumann algebra acting on the Fock space as follows
\begin{equation}\label{Rpm}
R(H) := \{\mathrm{w}(\xi): \xi\in H\}''\subseteq\cF(\cH)\ ,
\end{equation}
Let $H$ and $H_a$ be closed, real linear subspaces of $\cH$. Then we have
\begin{itemize}
\item[$(a)$] $R(\sum_a H_a) = \bigvee_a R(H_a)$,
\item[$(b)$] $R(\cap_a H_a) = \bigcap_a R(H_a)$,
\item[$(c)$] $R(H)' = R(H')$,
\end{itemize}
where $\bigvee$ denotes the von Neumann algebra generated (cf.~\cite{Ara63,LRT78,LMR16}).
We can now define the free scalar net of von Neumann algebras as
$$\cO\mapsto \cA(\cO):=R(\sH(\cO))\subseteq\cF(\cH).$$ It satisfies
\begin{itemize}
\item \textbf{Isotony:} if $\cO_1\subset\cO_2$ then $\cA(\cO_1)\subset \cA_U(\cO_2)$
\item \textbf{Locality:} if $\cO_1\subset\cO_2'$, then $\cA(\cO_1)\subset \cA(\cO_2)'$
\item \textbf{Poincar\'e covariance and positivity of the energy}: there exists a (anti-)unitary representation $\tilde U$ of $\cP_+$ on $\cF(\cH)$ such that $\tilde U(g)\Omega=\Omega$ and for $g\in\cP_+$, we have $ \tilde U(g) \cA(\cO)\tilde U(g)^*= \cA(gO)$. Furthermore the joint spectrum of the translation generators is contained in the forward lightcone $C$.
\item \textbf{Reeh--Schlieder:} $\overline{\cA(\cO)\Omega}=\cF(\cH)$ for any open $\cO\subset M$.
\end{itemize}

The free scalar net $\cA(\cO)$ satisfies the\textbf{ Bisognano--Wichmann property}, that claims
$$\tilde U(\exp(2\pi th_R))=\Delta_{\cA(W_R),\Omega}^{-it},\qquad t\in\RR$$ and the \textbf{modular reflection property}
\begin{equation}\label{eq:mrf}
\tilde U(\tau_{h_R})=J_{\cA(W_R),\Omega}.
\end{equation}

In even spacetime dimension the PCT Theorem ensures an antiunitary implementation of the space, time and charge reflection on the Hilbert space. Since the space and time reflection differs from $\tau_{h_R}$ by $\pi$-rotations, the modular reflection property is equivalent to the PCT-theorem, see  \cite{Mu01,BY01}. In odd spacetime dimension, the \eqref{eq:mrf} is considered as the PCT theorem (cf.~\cite{Mu10, BY01}).

\subsubsection{Nets of standard subspaces}\label{sect:stsub}
The mathematical structure of the real local subspaces is  recalled in this subsection.
We call a closed real subspace $\sH$ of the complex Hilbert space 
$\cH$ {cyclic} if $\sH+i\sH$ is dense in $\cH$, {separating} if $\sH\cap i\sH=\{0\}$, and \textit{standard} 
if it is cyclic and separating. We write $\Stand(\cH)$ for
the set of standard subspaces of $\cH$. The symplectic complement (or symplectic
orthogonal) of a real subspace $\sH$ is defined by the symplectic form 
$\Im \langle\cdot,\cdot\rangle$ on $\cH$ via 
\[ \sH'=\{\xi\in\cH: 
 \Im\langle\xi,\eta \rangle=0, \forall \eta \in \sH \}=(i\sH)^{\perp_{\Re\langle\cdot,\cdot\rangle}},\]
 where $\Re\langle\cdot\,,\cdot\rangle$ is the real part of the scalar product of $\cH$. 
Then $\sH$ is separating if and only if $\sH'$ is cyclic, hence $\sH$ is standard if and only if $\sH'$ is standard.
For a standard subspace $\sH$, one can
define the {\it Tomita operator} $S_\sH$ as the closed antilinear involution
\[S_\sH: \sH+i\sH \ni\xi+i\eta\longmapsto \xi-i\eta\in \sH + i \sH. \] 
The polar decomposition $S_\sH=J_\sH\Delta_\sH^{\frac12}$ defines an antiunitary involution $J_\sH$, called the \textit{modular conjugation} of $\sH$, and a positive self-adjoiont operator~$\Delta_\sH$, called the \textit{modular operator }of  $\sH$. 
Then we have that
\[  J_\sH\sH=\sH', \quad 
 \Delta^{it}_\sH\sH=\sH \qquad \mbox{ for every  } \quad 
t\in \R\] and
the following modular relation holds, namely
\begin{equation}\label{eq: TR}
J_{\sH}\Delta^{it}_{\sH}J_{\sH}=\Delta^{it}_{\sH}\qquad \mbox{ for every  } \quad t\in \R.
 \end{equation}
 
We remark that given a von Neumann algebra $\cA\subset\cB(\cH)$ with a cyclic and separating vector $\Omega\in\cH$  then one can define the standard subspace $\sH_\cA=\overline{\cA_{\text sa}\Omega}$, where $\cA_{\text sa}=\{a\in\cA: a=a^*\}$. By direct inspection one can see that $S_{\sH_\cA}=S_{\cA,\Omega}$, namely the standard subspace $\sH_{\cA}$ contains all the information about the Tomita operator of the von Neumann algebra $\cA$ with respect to the vector state $\Omega$.

One can completely reconstruct the standard subspaces from the related Tomita operators as \begin{equation}\label{eq:stsb}\sH=\Fix(S_\sH)=\ker(1-S_\sH)\end{equation}
(\cite[Thm.~3.4]{Lo08}).
This construction 
leads to a one-to-one correspondence between
{couples $(\Delta, J)$ where $J$ is an antiunitary operator  and $\Delta$ a selfadjoint positive operator, satisfying the modular relation:
  \begin{equation}\label{eq:commTT}
  J\Delta^{it}J=\Delta^{it},\qquad t\in\RR\end{equation}
  and standard subspaces $\sH=\Fix(J\Delta^\frac12)$,
see  \cite[Prop.~3.2]{Lo08}. Note that \eqref{eq:commTT} is equivalent to $J\Delta J=\Delta^{-1}$.

Standard subspaces have  a natural structure that adapts to the language of AQFT as one can see here 
\begin{itemize}
\item \textbf{Covariance property}:{{(\cite[Lem.~2.2]{Mo18})}.}
Let $\sH\subset\cH$ be a standard subspace  and $U\in\AU(\cH)$ 
be a {unitary or antiunitary} operator. 
Then $U\sH$ is also standard and 
$U\Delta_\sH U^*=\Delta_{U\sH}^{\epsilon(U)}$ and $UJ_\sH U^*=J_{U\sH}$, 
where $\epsilon(U) = 1$ if $U$ is unitary and 
$\epsilon(U) = -1$ if it is antiunitary. 

\item \textbf{Duality  property}(\cite[Prop.~2.13]{Lo08}):
Let $\sH\subset\cH$ be a standard subspace, and $(\Delta, J)$ the associated  operators, then the symplectic complement $\sH'$ is associated to the couple $(\Delta^{-1},J)$.
\item {\textbf{Inclusion properties }(\cite[Thms.~3.15, 3.17]{Lo08}, \cite[Thm.~3.2]{BGL02})}: 
Let $\sH\subset\cH$ be a standard subspace and $U(t)=e^{itP}$ 
be a unitary one-parameter group on $\cH$ with   generator $P$.
\begin{itemize}
\item[\rm(a)] 
If $\pm P > 0$ and $U(t)\sH\subset \sH$ for all $t\geq 0$, then 
\begin{equation}\label{eq:DeltaJ} \Delta_\sH^{-is/2\pi}U(t)\Delta_\sH^{is/2\pi} 
= U(e^{\pm s}t)\quad \mbox{ and } \quad 
J_{\sH}U(t)J_{\sH}=U(-t) \quad \mbox{ for all } \quad 
t,s \in \R.
\end{equation}
\item[\rm(b)] If $\Delta_\sH^{-is/2\pi}U(t)\Delta_\sH^{is/2\pi} 
= U(e^{\pm s}t)$ for $s,t \in \R$, then  the following are equivalent:
\begin{enumerate} 
\item[\rm(1)] $U(t)\sH \subset \sH$ for $ t\geq 0$;
\item[\rm(2)] $\pm P$ is positive.
\end{enumerate}
\end{itemize}
The first part of this result is the standard subspace version of the Borchers Theorem and  the second part is its converse.
\end{itemize}

We can now present an axiomatic framework for nets of standard subspaces on abstract wedges
\begin{definition}\label{def:ssn}Let $G$ be a connected Lie group with trivial center, such that $\cE(\fg)\neq \emptyset $. We fix an Euler element  $h\in\fg$, {then $\tau_h$ integrates to an involution on $G$}. Let $W_0=(h,\tau_h)$ be an Euler wedge and $C$ be an $\Ad^\epsilon(G_{\tau_h})$-invariant pointed convex cone in the Lie algebra $\fg$ of $G$.
Let $(U,\cH)$ be a unitary representation of $G$  and 
\begin{equation}
  \label{eq:net}
\sN \: \cW_+:=\cW(W_0) \to \Stand(\cH) 
\end{equation}
be a map  called a {\it net of standard subspaces on abstract wedges}. 
We consider the following properties: 
\begin{itemize}
\item[\rm(HK1)] {\bf Isotony:} $\sN(W_1) \subeq \sN(W_2)$ for $W_1 \leq W_2$, , $W_1,W_2 \in \cW_+$. 
\item[\rm(HK2)] {\bf Covariance:} $\sN(gW) = U(g)\sN(W)$ for 
$g \in {G}$, $W \in \cW_+$. 
\item[\rm(HK3)] {\bf Spectral condition:} 
$ C \subeq C_U := \{ x \in \g \: -i \partial U(x) \geq 0\}.$ 
We then say that $U$ is {\it $C$-positive}.   Note that $C_U$ is pointed if and only if $\ker(U)$ is discrete. 
\item[\rm(HK4)]{\bf Locality:} 
{If $W_1,W_2 \in \cW_{+}$ is such that $W_1\leq W_2'$, then $\sN(W_1) \subseteq \sN(W_2)'.$}
\item[\rm(HK5)] {\bf Bisognano--Wichmann (BW) property:}  for every $W=(h_W,\tau_W)\in \cW_+$
$$U(\exp(th_{W})) = \Delta_{\sN(W)}^{-it/2\pi},\qquad\forall t\in\RR,$$ 
\item[(HK6)] \textbf{Haag Duality}: $\sN(W')= \sN(W)'$ if $W, W' \in \cW_+$ 
\end{itemize}
If the representation $U$ extends antiunitarily to $G_{\tau_h}$ we can further 
require: 
\begin{itemize}
\item[(HK7)] \textbf{G-covariance:} 
{There exists an (anti-)unitary extension   of $U$ from $G$ to $G_{\tau_h}$ and an extension of the net $\sN(W)$ on $\cW:=\cW(W_0)$, such that the following condition is satisfied:}
$$ \sN(g. W) = U(g) \sN(W) \quad \mbox{ for } \quad 
g \in G_{\tau_h}, W\in\cW
$$
\item[(HK8)] \textbf{Modular reflection:} 
$U(\tau_{W})= J_{\sN(W)}$ for every $W=(h_W,\tau_W)\in\cW$.
\end{itemize}
\end{definition}

We remark that if $C=\{0\}$, there are no non-trivial wedge inclusions, then the isotony property trivializes and the locality property reduces to the condition $\sH(W')\subseteq\sH(W)'$ when $W,W'\in \cW_+$.

We further indicate the modular covariance property as the geometric action of the modular group: \begin{itemize}
\item[] \textbf{Wedge modular covariance}: for every couple of wedges $W_a,W_b\in\cW_+$,  
\begin{equation}\label{eq:mc}\Delta_{\sN(W_a)}^{-it}\sN(W_b)
=\sN(\Lambda_{W_a}(2\pi t).W_b)
\end{equation} for $W_a,W_b\in\cW_+$, $t \in \R$.
\end{itemize}
Note that the wedge modular covariance property does not require the Bisognano--Wichmann property, but it holds whenever the latter  is satisfied.

By identifying the abstract wedges with concrete wedge regions of  Minkowski space, the net of standard subspaces on wedges of the free scalar field given by \eqref{eq:ffn} satisfies all the previous assumptions.  

Lastly, we remark that an analogous set of axioms for nets of von Neumann algebras on abstract wedges that fits with second quantization can be provided  (\cite{MN21,MNO25}).

\subsubsection{The Brunetti--Guido--Longo (BGL) net}\label{sect:BGL}

The explicit construction of the free field we have seen in Sect.~\ref{sect:ff} for the free scalar field is highly non--canonical. It has been shown in \cite{LMR16} that such a construction is not feasible for infinite-spin representations of the Poincar\'e group as they lack localized states in bounded regions. However,  wedge subspaces are determined by the wedge symmetries. 
Specifically, given an  (anti-)unitary representation of the proper Poincar\'e group $\cP_+$, for any wedge $W\subset\RR^{1,d}$, the wedge subspace is reconstructed as follows:
The operators
$$\Delta_W^{it}:=U(-2\pi h_W t),\qquad J_W:=U(\tau_{h_{W}})$$ 
satisfy the commutation relation \eqref{eq:commTT}, hence $\Delta_W=\exp(2\pi i \partial U(h_W))$ and $ J_W=U(\tau_{h_{W}})$  define a the wedge standard subspace $$\sH(W)=\Fix(J_W\Delta_W^{\frac12})$$ as we have seen in the previous section. Based on these remarks R.~Brunetti, D.~Guido and R.~Longo in \cite{BGL02} provided a construction of the one-particle net of standard subspaces of the free Poincar\'e and conformal covariant theories on Minkowski space, Lorentz covariant theory on de Sitter space and M\"obius covariant theory on the chiral circle, on the corresponding index sets of wedge regions (see also \cite{GLW,Lo08}). 
This construction has been generalized to our abstract set of wedges for a $\Z_2$-graded Lie group supporting Euler elements in \cite{MN21}. The anti-de Sitter space can be included in the picture and it is treated in \cite{NOO21,NO23}, see also \cite{Re00}.

Let $G$ be a connected Lie group with trivial center, such that $\cE(\fg)\neq \emptyset$. {We fix an Euler element  $h\in\fg$ and  $\tau_h$ integrates to an involution on $G$.}
Let $(U,\cH)$ be an (anti-)unitary representation of $G_{\tau_h}$,   $W  = (h, \tau_h) \in \cG_E$ be an Euler wedge, and consider the couple of operators
\begin{equation}
  \label{eq:bgl}
J_{\sH_U(W)} = 
U(\tau_h) \quad \mbox{ and } \quad \Delta_{\sH_U(W)} = e^{2\pi i \partial U(h)}.
\end{equation}
They satisfy \eqref{eq:commTT} and we can define $\sH_U(W):=\Fix(J_{\sH_U(W)}\Delta_{\sH_U(W)}^{\frac12})$. 
The map
  \[ \sH_U^{\rm BGL} \: \cG_E(G_{\tau_h}) \ni W\to \sH_U^{\rm BGL}(W)\in\Stand(\cH)\]
  is the so-called BGL net associated to $U$ and satisfies all the assumptions in Definition \ref{def:ssn}.

  \begin{theorem} \mlabel{thm:BGL}\cite[Thm~4.12, Prop.~4.16]{MN21}
{    Let $(U,\cH)$ be an (anti-)unitary C-positive representation of $G_{\tau_h}$ where $h\in \cE(\fg)$. }
    Then the BGL net 
\[ \sH_U^{\rm BGL} \: \cG(G_{\tau_h}) \to \Stand(\cH) \] 
 satisfies (HK1)-(HK8) in Definition \ref{def:ssn}.
\end{theorem}

Theorem \ref{thm:BGL} can be stated in a more general framework of abstract non-Euler wedges as in \cite{MN21}.  
Given a $\Z_2$-graded Lie group $G_\sigma=G\rtimes \{e,\sigma\}$ with $\sigma$  an involution on $G$, then one can define more general notion for the set of abstract wedges without referring to Euler elements as follows:
\begin{equation}\label{eq:genwed}
\cG:=\{(h,\tau)\in \fg\times  G\sigma: \Ad (\tau)(h)=h, \tau^2=e\}
\end{equation} 
This set satisfies analogous properties to the set of Euler wedges and the BGL net is still well defined on $\cG$ given an (anti-)unitary representation of $G_\sigma$ \cite[Sect.~2]{MN21}. In Section \ref{subsec:3.1} we will explain why Euler elements, hence Euler wedges, are the natural choices for an index set of abstract wedges.

We have introduced an axiomatic framework for nets of standard subspaces on wedges. This picture includes the well known models from AQFT and produces new examples. Given any three graded Lie algebra $\fg$ and $h\in\cE(\fg)$, the group $G=\Inn(\fg)$ supports unitary representations that extend (anti-)unitarily to $G_{\tau_h}$ (up to properly coupling representations). Consider $\cG_E(G_{\tau_h})$ the set of Euler wedges and  $\Ad^\eps (G_{\tau_h})$-invariant pointed cone $C$ in $\fg$ defining the inclusion rules of Euler wedges.
 Then one can construct the BGL-net of standard subspaces on abstract wedge regions. We stress that if the cone $C=\{0\}$, namely there are no non-trivial wedge inclusion, then isotony property is trivially satisfied. Hermitian simple Lie algebras have a non-trivial positive cone and their unitary positive energy representations are studied in \cite{MdR07}. Then if $U$ is a positive energy representation satisfying the spectral condition (HK3) with respect to $C\neq\{0\}$, the isotony property of the BGL-net follows by the inclusion property of standard subspaces, in particular by the converse of the Borchers theorem (cf. Sect.~\ref{sect:stsub}, \cite[Thm.~4.12]{MN21}). Once a net of standard subspaces on abstract wedge regions is constructed the  second quantization net of von Neumann algebras is canonically provided by relations \eqref{eq:weyl} and \eqref{eq:vs}.  If a Lie group  $G$ has a non-trivial center, the locality property has to be reformulated in terms of twisted central complements. This has been investigated in \cite[Sect.~4.2.2]{MN21}. In this case a proper second quantization procedure respecting a twisted locality conditions has still to be established. Recent results on second quantizations are contained in \cite{CSL23}.

We have seen that there is a correspondence between abstract wedges and concrete wedges on causal manifolds. It is possible to provide a construction for a net of standard subspaces directly on the causal manifold extending the construction presented for the free field on Minkowski space. 
  
Given a unitary representation $(U,\cH)$ of a connected a Lie group~$G$ 
and a homogeneous space $M = G/H$, we are interested in
families $(\sH(\cO))_{\cO \subeq M}$ of closed real subspaces of $\cH$, 
indexed by open subsets $\cO \subeq M$,  called {\it nets of real subspaces on $M$}. 

This setting is depicted by the following set of axioms:
\begin{itemize}
\item[(Iso)] {\bf Isotony:} $\cO_1 \subeq \cO_2$ 
implies $\sH(\cO_1) \subeq \sH(\cO_2)$ 
\item[(Cov)] {\bf Covariance:} $U(g) \sH(\cO) = \sH(g\cO)$ for $g \in G$. 
\item[(RS)] {\bf Reeh--Schlieder property:} 
$\sH(\cO)$ is cyclic  if $\cO \not=\eset$. 
\item[(BW)] {\bf Bisognano--Wichmann property:}
{There exists an open subset $W \subeq M$ (called a {\it wedge region}),  
such that $\sH(W)$ is standard
and $\Delta_{\sH(W)} = e^{2\pi i \partial U(h)}$
  for some $h \in \g$.}
\end{itemize}
{ We remark that locality is not present among the list, since it is to be investigated in \cite{MN25}.}
Nets satisfying (Iso), (Cov), (RS), (BW) 
 have been constructed on non-compactly causal symmetric spaces in \cite{FNO23} and on compactly causal spaces in \cite{NO23}, using the language of distribution vectors.
Given the $G$-orbit of a wedge region, one can define the \textit{maximal net }on bounded regions by intersection of wedge subspaces $\sH^{\text{max}}(\cO)=\bigcap\{\sH(gW): \cO\subset gW, g\in G \}$, this has been studied in \cite[Sect.~2.2.4]{MN24}.

\section{Euler elements, Bisognano--Wichmann property and algebra types}\label{sect:3}
\subsection{Non-modular covariant nets of standard subspaces}
\label{sect:modEul}

A long standing question in AQFT concerns the necessity to assume the Bisognano--Wichmann property. In particular, given an isotonous, covariant, local, (positive energy) net of von Neumann algebras, does the Bisognano--Wichmann property hold? This is always the case if the local algebra net is generated by Wightman fields as it is proved in \cite{BW75}. It holds  in conformal theories  (\cite{Lo08,BGL93}) and in massive theories \cite{Mu01}. One can give an algebraic sufficient condition on the covariant representation ensuring the Bisognano--Wichmann property for one-particle nets that can be applied to prove the Bisognano-Wichmann property in some possibly interacting theories  \cite{Mo18,DM20}. When infinite particle degeneracy is present many counterexamples to Bisognano-Wichmann property  arise \cite{Mo18,MT19}. On the other hand in these counterexamples  the wedge modular groups implement a covariant representation of the symmetry group. 
 In \cite{Y94}, Yngvason provided a translation covariant von Neumann algebra net  where there is no geometric modular action of the wedge algebra.  On the other hand these nets are not expected to be Lorentz covariant. {Further developments about  the relation between Modular Action, Wedge Duality and Lorentz Covariance for generalized free fields are in \cite{GY00}}.

Since the modular theory of a net of von Neumann algebras in the vacuum representation is encoded in the standard subspace structure of wedge algebras, it is a natural problem to provide counterexamples to the Bisognano--Wichmann property in the net of standard subspaces context without the modular covariance property  of  wedge subspaces. This problem has been recently faced in \cite{MN22} providing a family of non-local standard subspace nets with no geometric action of the wedge modular groups. We summarize below the main steps of the construction.

{Let $G=\Inn(\fg)$ where $\fg$ is simple, let $G_\tau=G\rtimes\{1,\tau\}$, where $\tau$ is an Euler involution, and let $U$ be an {(anti-)unitary} representation of $G_\tau$ on a Hilbert space $\cH$.}
 Assume the existence of a subgroup $H_\tau=H\rtimes\{e,\tau\} \subeq G_\tau$, 
two Euler wedges $W_1 = (h_1, \tau_1)\in \cG_E(H_\tau)$ and $ W_2 = (h_2, \tau_2) \in \cG_E(G_\tau)$. 
Note that $h_1$ is not necessarily in $\cE(\fg).$
Assume that the stabilizer $H_{W_1}$ of $W_1$ in $H$ fixes $W_2$, 
in particular $[h_1, h_2] = 0$. We denote the orbit $H.W_1$ with $\cW_+(H,W_1)$.

Consider the BGL construction with respect to $U$, $G_\tau$ and $\cG_E(G_\tau)$. The standard subspace $\sN_2 = \sN_U(W_2)$ is defined by
\[ J_{\sN_2} = U(\tau_2) \quad \mbox{ and  }\quad 
\Delta_{\sN_2} = e^{2 \pi i \cdot \partial U(h_2)}.\] 
Since the stabilizer $H_{W_1}$ fixes $W_2$ and the unitary group $U(H_{W_1})$ fixes $\sN_2$, and
we can define an  $H$-equivariant map 
\[ \sN \: \cG_E(H_\tau) \supeq \cW_+ := 
\cW_+(H,W_1) := H.W_1  \to \Stand(\cH), \quad 
g.W_1 \mapsto U(g)\sN_2\] 
which is uniquely determined by 
\begin{equation}
  \label{eq:incond}
 \sN(W_1) = \sN_2.
\end{equation}

The net $\sN$ is not necessarily modular covariant. 
The following is an equivalent condition for the {modular covariance} property \eqref{eq:mc}  (cf.~\cite[Lem~3.1]{MN22}).
The net $\sN$ on $\cW_+(H,W_1)$ satisfies  modular covariance 
if and only if, for all ${g \in H},  t \in \R$, the operator 
\[ U(g) U(\exp t (h_1-h_2)) U(g)^{-1} \] 
fixes the standard subspace $\sN_2$, i.e., 
\begin{equation}
  \label{eq:covcond1}
 g \exp(t (h_1-h_2))g^{-1} \in G_{\sN_2} \quad \mbox{ for }\quad 
g \in H, t \in \R. 
\end{equation}
If $\ker(U)$ is discrete, then  \eqref{eq:covcond1} is violated if \begin{equation}\label{eq:eq:covcond2}
[h_1-h_2,\fh]\not\subseteq \ker(\ad\,h_2),\end{equation}
where $\fh=\Lie(H)$, see \cite[Rem.~3.2]{MN22}.

We now indicate how to use this prescription to construct  a {non-modular covariant net on the two-dimensional de Sitter space.} 
The crucial step here is to choose $h \in \mathfrak{g} $ as a \emph{non-symmetric} Euler element in $ \mathfrak{g} $, and to verify the failure of condition~\eqref{eq:covcond1} by testing condition~\eqref{eq:eq:covcond2}.
 We refer to \cite{MN22} for the general discussion and further examples as non-modular covariant nets on Minkowksi space and the representation theoretic results.

In the context described above, we fix $\fg$ to be a simple non-compact real Lie algebra. 
Consider a {\bf non-symmetric} Euler element $h_2\in\fg$ 
and the associated wedge $W_2=(h_2,\tau_2) \in \cG_E(G_{\tau_2})$.
By \cite[Lem.~2.15]{MN22} there exists a 
$\fgl_2$-subalgebra $\fb\subeq \g$ containing $h_2$  
such that  $\fh := [\fb,\fb] \cong \fsl_2(\R)$ 
satisfies  $[\fh,h_2]\neq 0$. We consider 
$H := \Inn_\g(\fh) \cong \SL_2(\RR)$ and $H_{\tau_2}=H\rtimes\{1,\tau_2\}$ (see \cite[Lem~2.15(d)]{MN22}).
{The Euler element $h_2 \in \fb$ has  a central component~$h_c$, so that
  \begin{equation}
    \label{eq:hsplit1}
    h_2 = h_c - h_1 \quad \mbox{ with } \quad
h_c \in \fz(\fb) \quad \mbox{ and } \quad
h_1 \in {\cE(\fh).}
  \end{equation}
Choosing the isomorphims $\fh\to\fsl_2(\RR)$ suitably, 
 we find
\begin{equation}
  \label{eq:hsplit2}
  h_1 =\frac12\left(\begin{array}{cc}0&1\\1&0
\end{array}\right).
\end{equation}}

Let $V$ be the restriction $U|_{H_{\tau_2}}$. 
The group  $H\cong \SL_2(\RR)$ is the double covering of $\PSL_2(\R) 
\cong \cL_+^\up$, the Lorentz group of $\RR^{1,2}$. So it 
acts on de Sitter space $\dS^2$ 
through the covering map 
$\Lambda: \SL_2(\RR)\ni g\mapsto \Lambda(g)\in \PSL_2(\RR)$. We thus obtain a net $\sH^{\dS}$ on de Sitter spacetime as follows:
Let
     \[ {\cW^{\rm dS}}\ni W^{\rm dS}\mapsto \sH^{\dS}(W^{\rm dS})\subset\cH_U \]
such that
\begin{equation}
  \label{eq:exgen1}
  \sH^{\dS}(\Lambda(g)W_{R}^{\rm dS}):=V(g)\sN_U(W_2)\quad 
  \mbox{ for  } \quad g \in H \cong \SL_2(\R),
\end{equation}
where $\sN_U$ is the BGL net defined by $U$  \cite[Thm.~3.4]{MN22}. Then the net $\sH^{\dS}$ is Lorentz covariant and does not satisfy modular covariance. Indeed $h_2-h_1=h_c-2h_1 $ and
$$[h_2-h_1,\fh]=[h_1,\fh]\not\subseteq \ker(\ad\,h_1).$$ In particular \eqref{eq:eq:covcond2} is violated since $\ad h_2|_{\fh}=-\ad h_1$ and $\sN$ is not modular covariant.
One can also see this by observing that the Lie algebra generated by $gh_2g^{-1}$ with $g\in\SL_2(\RR)$ is $\fgl_2(\RR)$ so that $\Delta^{it}_{\sH^{\dS}(g W^{\dS}_{R})}=U(\exp(-2\pi gh_2g^{-1} t))$ can not generate a representation of $\tilde\SL_2(\RR)$.  Details in \cite[Sect.~3.2.1 and 3.2.2]{MN22}.
We can further comment that  counterexamples satisfying Lorentz covariance,  locality and violating modular covariance are not covered by these examples and are still missing.

\subsection{{The Euler Element Theorem}}
\mlabel{subsec:3.1}

The following theorem is an important step in the understanding the key role of Euler elements in these constructions.
It addresses the question: Are Euler elements necessary for defining fundamental localization regions in one-particle and von Neumann algebra nets? Given a unitary representation of a Lie group and a standard subspace $\sV$ satisfying the Bisognano--Wichmann property with respect to a  generic one-parameter group in $G$,  the generator of the latter in $U$ is an Euler element if a regularity property is satisfied.  

\begin{thm} \mlabel{thm:2.1} {\rm(Euler Element Theorem)} \cite[Thm.~3.1]{MN24}
  Let $G$   be a connected finite-dimensional Lie group with 
  Lie algebra $\g$ and $h \in \g$. 
  Let $(U,\cH)$ be a unitary 
  representation of $G$ for which $$\ker(\dd U)\cap [h,\fg]=\{0\}.$$ 
  Suppose that $\sV$ is a standard subspace
    and $N \subeq G$ an identity neighborhood such that the following hold:
  \begin{itemize}
  \item[\rm(a)] \textbf{Bisognano--Wichmann property:} $U(\exp(t h)) = \Delta_\sV^{-it/2\pi}$ for $t \in \R$.
  \item[\rm(b)] \textbf{Regularity property}: $\sV_N := \bigcap_{g \in N} U(g)\sV$ is cyclic.
  \end{itemize}
  Then $h$ is an Euler element (or central) and the conjugation $J_\sV$ satisfies
  \begin{equation}
    \label{eq:J-rel}
 J_\sV U(\exp x) J_\sV = U(\exp \tau_h(x)) \quad \mbox{ for } \quad
 \tau_h = e^{\pi i \ad h}, x \in \g.
  \end{equation}
\end{thm}
{The theorem implies that, under the above assumptions,  if $\g$ is a compact Lie algebra, then 
    $h \in \g$ is central,} so that ${\tau_h = \id_\g}$.
 Therefore, in this case, a standard subspace $\sV$ associated
    to a pair $(h,\tau) \in \cG_E$  by the BGL construction  
    can only satisfy the regularity condition in Theorem~\ref{thm:2.1}(b)
  if $\sV$ is $U(G)$-invariant.

It is a consequence of Theorem~\ref{thm:2.1} 
     that $\tau_h$ integrates to an
    involutive automorphism {$\tau_h^{U(G)}$ on the group $U(G) \cong G/\ker(U)$
    that is uniquely determined by
\[    \tau_h^{U(G)}(\exp x) = \exp(\tau_h(x)) \quad \mbox{ for }\quad x \in\g.\]
In particular the theorem ensures the extendibility of $U$ to an (anti-)unitary representation of $G_{\tau_h^G}$ on the same Hilbert space if $U$ is a faithful representation of $G$ satisfying Bisognano--Wichmann and regularity properties.}

 The regularity property is very natural in the setting of an analysis of state localization in AQFT, see for instance \cite{BGL02, GL95, LMR16}. 
 Firstly note that the condition (b) is equivalent to the existence of  a cyclic subspace
    $\sK\subset \sH$ such that $U(g)\sK\subset \sV$ for every $g\in N$. This corresponds to a localization property of  standard subspace nets that are finer than wedges. To provide examples we firstly refer to theories on Minkowski spacetime.
We call an open subset $\cO\subset\RR^{1,d}$ {\it spacelike} if $x_0^2 < \bx^2$ holds for all $(x_0,\bx) \in \cO$. Here we will consider spacelike cones that are also pointed, convex and causally complete cones, namely $\cC=\cC''$.
Given a   spacelike cone $\cC$ with apex in $a\in\RR^{1,d}$ then $\cC=\bigcap\{ gW_R+a: g\in\cL_+^\up, gW_R\supset\cC-a\}$ and for every wedge $W$ there always exists a spacelike cone $\cC\subset gW$, for every $g$ in some small enough neighbourhood of the identity $N\subset \cP_+^\up$.
 Let $\sH$ be a net of real subspaces on open regions satisfying isotony, covariance,  cyclicity on wedge regions, then  the assumption (a) and (b) in the theorem correspond respectively to the Bisognano--Wichmann property and the cyclicity for  $\sH(\cC)$. The latter property is satisfied by the free one-particle nets on Minkowksi space when $\cC$ is a spacelike cone, see  \cite[Sect.~4]{BGL02}.   One can make an analogue remark for de Sitter spacetime. One can consider the intersection of pointed convex spacelike cones with apex in the origin in $\RR^{1,d}$ with  de Sitter space $\rm{dS}^d$ and $\cC\cap\rm{dS}^d$ are bounded regions obtained by wedge intersection.

Another remark, in view of \cite{DM20} there exists no net of standard subspaces on spacelike cone regions, satisfying the Bisognano-Wichmann property for wedge subspaces, undergoing an irreducible massless finite (non-zero) helicity $\cP_+^\up$-representation $U$. By contradiction, if there exists such a net, then Theorem \ref{thm:2.1} applies and $U$ should extend to  a representation of $\cP_+$ on the same Hilbert space. This is never the case for finite non-zero helicity representation, cf. \cite[Thm.~9.10]{Var85}.

The Euler Element Theorem can be applied to nets of standard subspaces on causal manifolds, see \cite[Thm.~3.4]{MN24}. With notations as in Sect.~\ref{sect:BGL},  let $(U,\cH)$ be a unitary representation of the
    connected Lie group $G$ with $\ker(U)$ discrete.
    If $(\sH(\cO))_{\cO\subeq M}$ is a net of real subspaces on
   {(the open subsets of)} a $G$-manifold $M$ that satisfies
        {\rm(Iso), (Cov), (RS)} and {\rm(BW)}, then Theorem \ref{thm:2.1} applies to $\sH(W)$, the Lie algebra element $h$  is an Euler element, and
    the conjugation $J := J_{\sH(W)}$ satisfies \eqref{eq:J-rel}.

The theorem  applies also when von Neumann algebras are taken into account. 
 Let $G$, $h$ and $(U,\cH)$ be as in Theorem \ref{thm:2.1}. Let $\Omega$ be a unit vector  and   $\cM \subeq B(\cH)$ be a von Neumann algebra
  for which $\Omega$ is cyclic and separating.
Assume that $\Omega$ is fixed by $U(G)$ and that the Bisognano--Wichmann property holds as
  $U(\exp({-2\pi t h})) = \Delta_{\cM,\Omega}^{it}$.
Assume that for some $e$-neighborhood $N \subeq G$, the vector
  $\Omega$ is still cyclic for the
  von Neumann algebra 
  \[ \cM_N := \bigcap_{g \in N} \cM_g, \quad \mbox{ where } \quad
    \cM_g= U(g)\cM U(g)^{-1}.\]
  Then $h$ is an Euler element and the modular conjugation $J = J_{\cM,\Omega}$
of the pair $(\cM,\Omega)$ satisfies
  \[ J U(\exp x) J = U(\exp \tau_h(x)) \quad \mbox{ for } \quad
  \tau_h = e^{\pi i \ad h}.\]
This follows by applying Theorem \ref{thm:2.1} to the standard subspace $\sH=\overline{\cM_{sa}\Omega}$, see \cite[Thm.~3.7]{MN24}.

Borchers and Buchholz in \cite[Theorem~6.2]{BB99} proved that given a local, Lorentz covariant, weakly additive net of von Neumann algebras on bounded regions of de Sitter space,  assuming that vacuum state looks for any observer in a wedge moving through the corresponding boost-time evolution like equilibrium state with some apriori arbitrary temperature, namely a KMS-state, then the geodesic temperature of vacuum state has to be equal to the Gibbons--Hawking temperature. This corresponds to the fact that the boost generator - properly parametrized - is an Euler element. A first geometric approach to the \cite{BB99} result is in \cite{Str08}. The Euler Element Theorem  gives a standard subspace reformulation of this picture \cite[Cor.~5.16, 3.3 and 3.5]{MN24}.

\subsection{Regularity and Localizability}\mlabel{sec:loc-reg}

Given a connected  Lie group $G$, such that $h\in\cE(\fg)$,  and  an (anti-)unitary representation  $U$  of $G_{\tau_h}$, one can construct a net of standard subspaces on wedges  as in Sect.~\ref{sect:BGL}. The Bisognano–Wichmann property holds, which naturally raises the question of the extent to which the regularity property (b) from Theorem \ref{thm:2.1} is also satisfied. 
This is a representation theoretic question investigated in \cite{MN24} that in part uses the results of \cite{FNO23}. 

In this setting we reformulate the property as follows. Assuming that $h$ is an Euler element and $(U,\cH)$ an (anti-)unitary representation 
  {of $G_{\tau_h}$}. { Let $\sV$ be the standard subspace associated to the wedge $(h, \tau_h)$ by the BGL--construction}. We shall say that $U$ is {\it regular with respect to $h$}, or
{\it $h$-regular}, if there exists an 
$e$-neighborhood $N \subeq G$ such that
$\sV_N = \bigcap_{g \in N} U(g)\sV$ is cyclic.
{Replacing $N$ by its interior, we may always assume that $N$ is
  open.}

This property passes to direct sums, direct integrals and sub-representations as follows \cite[Lem.~4.4]{MN24}:
  \begin{itemize}
  \item[\rm(a)] If $U = U_1 \oplus U_2$ is a direct sum, then
    $U$ is $h$-regular if and only if $U_1$ and $U_2$ are $h$-regular.
  \item[\rm(b)] If $U$ is $h$-regular, then every subrepresentation is
    $h$-regular. 
    \item[\rm(c)] A  direct integral representation $U = \int_X^\oplus U_m \, d\mu(m)$
  is regular if and only if there exists an
    $e$-neighborhood $N \subeq G$ such that, for $\mu$-almost every $m \in X$,
    the subspace $\sV_{m,N}$ is cyclic.
  \end{itemize}

We  now   present the state of the art on the analysis of this property.

\subsubsection{Regularity via generating positive cone}
In conformal Algebraic Quantum Field Theory, the geometry of the model ensures that one can shrink wedge regions:  conformal symmetries in Minkowski space can transform a wedge region $W$ into a doublecone region $O\Subset W$ and both regions are associated 1-1 to Euler elements; for  the M\"obius wedge geometry, intervals are associated to Euler wedges, and interval inclusions are associated to inclusions of Euler wedges \cite{BGL93,BGL02,MN21}. In particular, provided the BGL construction for these models, there exists an open neighborhood of the identity $N\subset G$ such that given a $\bigcap_{g\in N}gW_1 $ contains a wedge region supporting a cyclic subspace, hence regularity for the BGL net holds.  

{Firstly, we situate the general framework within the Lie theory context, where the  regularity property stems from the existence of a generating positive cone}. Theorem \cite[Thm.~4.9]{MN24}  claims that, if $(U,\cH)$ is an (anti-)unitary representation of $G_{\tau_h}$
  for which the cones
  \begin{equation}\label{eq:gpc}    C_\pm := \pm C_U \cap \g_{\pm 1}(h)\end{equation}
  are linearly generating  $\g_{\pm 1}(h)$, then $(U,\cH)$ is regular.

As examples one can consider the M\"obius group. Here, the condition \eqref{eq:gpc} holds for
{positive energy representations of}
  the M\"obius group. Up to sign, the
  only pointed, generating closed convex $\Ad$-invariant cone has been introduced in \eqref{eq:Csl}.
Given a positive energy representation $U$ of $\Mob$, then $C_U=C$. Considering the Euler element $\displaystyle{h = \frac{1}{2}\pmat{1 & 0 \\ 0 & -1}}$, we have
\[ C_\pm := \pm C \cap \g_{\pm 1}(h), \quad
  C_+ = \R_+ \pmat{0& 1\\ 0 & 0}, \quad 
  C_- = \R_+ \pmat{0& 0\\1 & 0},\]
and the half lines $C_\pm$ generate $\g_{\pm 1}(h)$.  
This property is not satisfied by positive energy representations of the  Poincar\'e group 
on $\R^{1,3}$. We will see in the next section how to deal with this case.

Another example is the the following. Let $\fg=\so(2,d+1)$ the Lie algebra of the conformal group of the Minkowski spacetime $\RR^{1,d}$,  consider the Euler element $h\in\fg$ generating the dilations $\exp(ht)x=e^tx$ where $x\in \RR^{1,d}$, then $\fg_1\simeq\RR^{1,d}$ and the positive cone $C_+=\{x\in\RR^{1,d}:x^2\geq0,x_0\geq0\}$ is generating in $\fg_1$ (and analogously for $\fg_{-1}$). {The generating property holds for all hermitian simple real Lie algebras, see \cite[Lem.~3.2]{NOO21}.}
 
\subsubsection{Regularity in a semidirect product} 

As a second step, it looks natural to study the regularity property for (anti-)unitary representations of semidirect products where one can split the problem checking two conditions on the two defining subgroups. This picture now contains the Poincar\'e group that is not simple. 
{Furthermore any simply connected Lie group $G$,  is a semidirect product $G \cong N \rtimes S$, where $S$ is semisimple and $N$ is the solvable radical (Levi decomposition, see e.g.~\cite[Thm. 11.1.19]{HN12}). 

We consider a general (connected) semidirect product Lie group $G = R \rtimes L$, let $ h\in \cE(\fg)$, and $(U,\cH)$ be an (anti-)unitary representation of $G_{\tau_h}$. }Theorem 4.11 in \cite{MN24} proves that $(U,\cH)$ is regular under the following two conditions
  \begin{itemize}
 \item[(a)] the cones $C_\pm := \pm C_U \cap \fr_{\pm 1}(h)$ generate
    $\fr_{\pm 1}(h)$ where $\fr=\Lie (R)$. 
  \item[(b)] $(U\res_{L_{\tau_h}},\cH)$ is $h$-regular.
   \end{itemize}

{This theorem applies to (anti-)unitary positive energy representations of the Poincar\'e group $G_{\tau_{W_R}}$. In this context let $R=\RR^{1,d}$, $L=\cL_+^\uparrow$ and $h=h_{W_R}\in\cE(\fg)$, one can check the generating property of the one-dimensional cones $C_\pm$ in the eigenspaces
    $\fr_{\pm 1}(h_{W_R}) $ that corresponds to the lightrays $ \R (\be_0 \pm \be_1)\fr=\fr_{\pm1}$  and the regularity property for the restriction of the representation
  to the {identity component $\cL_+^\up$ of the Lorentz group.}
  The first property follows from the spectral condition, namely positivity of the energy of the the Poincar\'e representation. The second one holds for every representation of the Lorentz group, as one can see from Theorem 4.25 in \cite{MN24}. 
  
To study the regularity property for reductive Lie groups $G$ one can introduce the localizability property. Given an (anti-)unitary representation $(U,\cH)$ of the Lie group $G_{\tau_h}$, this property ensures the existence of a concrete cyclic real subspace of the representation Hilbert space  localized on an open region of a causal symmetric manifold $M$ lying in the intersection of wedge subspaces (wedge regions as described in Sect.~\ref{sect:chs}). 
In \cite{FNO23}, the authors explicitly constructs local subspaces (in analogy with the the free field case) and prove the Bisognano--Wichmann property for standard subspaces of wedge regions and the cyclicity of the subspaces associated to open regions that are intersections of wedge regions. This result is recalled to conclude regularity in \cite{MN24}. 
  Results in \cite[Thm.~4.23, Cor~4.24]{MN24} guarantee localizability
  for representations of a connected reductive Lie group.  Then Theorem 4.11 in \cite{MN24} applies to the cases where $G=R\rtimes L$, $U|_R$ satisfies (a) and $L$ is reductive.
  
  Currently, it is unknown whether all (anti-)unitary representations of Lie groups of the form $G_{\tau_h}$ are regular. The preceding discussion suggests that resolving this question requires a more detailed analysis of the case of solvable groups.

\subsection{The type of wedge algebras}
\label{sect:vna}
Von Neumann algebras  can be classified in terms of the geometry of their projection and {the spectrum of the modular action }\cite{Ta02,Su87, Co73}. {Algebras involved in an AQFT are proven to be of type III$_1$ in many cases, and they are expected to always be of this type when nontrivial (see for instance \cite{Dr77, Lo82, Fr85, BDF87, BB99, FiGu}).} Recently, type II$_1$ algebras have become increasingly relevant in certain constructions within the AQFT context, see \cite{CLPW23,FJLRW25}.

{We prove that our generalized AQFT setting  also respects this type III$_1$ expectation.} 
{Let $G$ be a connected Lie group with
    Lie algebra $\g$. We say that $h$ is \textit{anti-elliptic} if $\fn_h + \R h = \g$ where $\fn_h \trile \g$ be the smallest ideal of $\g$ such that the image of $h$ in
the quotient Lie algebra $\g/\fn_h$  is elliptic.}  
  Let $(U,\cH)$ be a unitary 
  representation of $G$ with discrete kernel,
  $\cN\subset \cM \subeq B(\cH)$ an inclusion of von Neumann algebras, 
  and $\Omega\in\cH$ a unit vector which is
  cyclic and separating for~$\cN$ and $\cM$. 
Assume that the following properties hold: 
  \begin{itemize}
  \item[\rm(Mod)] $e^{2\pi i \partial U(h)} = \Delta_{\cM,\Omega}$,  and 
  \item[\rm(Reg')] $\{ g \in G \: \Ad(U(g)) \cN \subeq \cM\}$
    is an $e$-neighborhood in $G$.
  \end{itemize}
  Then the following assertions hold:
  \begin{itemize}
  \item[\rm(a)]   $h$ is an Euler element. 
  \item[\rm(b)]  The conjugation $J := J_{\cM,\Omega}$ satisfies
  \begin{equation}\label{eq:jcov} J U(\exp x) J
    = U(\exp \tau_h(x)) \quad \mbox{ for } \quad
    \tau_h = e^{\pi i \ad h}, x \in \g.\end{equation}
\item[\rm(c)] $\cH^G = \ker(\partial U(h))$.   
\item[\rm(d)] The restriction of $i\partial U(h)$ to the
  orthogonal complement of the subspace 
  $\cH^{N_h}$ of fixed vectors of the {codimension-one} normal subgroup~$N_h:=\langle\exp(\fn_h)\rangle\trile G$,
  has absolutely continuous spectrum. 
  \end{itemize}
 If, in addition, $\cH^G = \C \Omega\not=\cH$, 
  then $\cM$ is a factor (i.e.~$\cM\cap\cM'=\bC\textbf{1}$) of type {\rm III}$_1$. This theorem is \cite[Thm.~5.15]{MN24}.
  
Note that if $h \in \cE(\g)$ is anti-elliptic, then $\fg_0\subseteq\RR h+[\fg_1,\fg_{-1}]$.  Anti-elliptic elements are called essential elements in \cite{Str08}.  
One can  apply this setting when it is considered a 
 unitary representation $U$ on an Hilbert space $\cH$
  of the {connected} Lie group
  $G$ with discrete kernel, a vector $\Omega\in\cH$ that is the unique unit vector fixed by $U(G)$ and
  $\cM \subeq B(\cH)$ be a von Neumann algebra
  for which $\Omega$ is cyclic and separating. Assume that $h\in\fg$ is anti-elliptic and the following properties.
  \begin{itemize}
  \item[\rm(Mod)] \textbf{Modularity:} There exists an  element $h \in \g$ for which
  $e^{2\pi i \partial U(h)} = \Delta_{\cM,\Omega}$. 
\item[\rm(Reg)] \textbf{Regularity:} For some $e$-neighborhood $N \subeq G$, the vector
  $\Omega$ is still cyclic (and obviously separating) for the
  von Neumann algebra 
  \[ \cM_N := \bigcap_{g \in N} \cM_g, \quad \mbox{ where } \quad
    \cM_g= U(g)\cM U(g)^{-1}.\]
\end{itemize}
Then $h$ is an Euler element and $\cM$ is a type III$_1$ factor. 
This picture includes the cases of nets of von Neumann algebras on Euler wedges when, under the previous assumptions, the map  $$\cW_+\ni W\mapsto \cM(gW)=U(g)\cM(W_h)U(g)^{-1}\subset\cB(\cH),\qquad g\in G $$ where $\cM(W_h)=\cM$ and $\,W_h=(h,\tau_h),$ defines a $G$-covariant isotonous net of von Neumann algebras.

In the case $\Omega$ is cyclic and separating but the set of $U(G)$-fixed points is not one dimensional, then one can consider the von Neumann algebra $\displaystyle{\cA := \big(\bigcup_{g \in G} \cM_g)'' \subeq B(\cH)}$
   generated by all algebras $\cM_g=U(g)\cM U(g)^{-1}$. If $\cM'=\cM_{g_0}$ for some $g_0\in G$ then we have direct integral decompositions
  \[ \cM=\int_X^\oplus\cM_x\, d\mu(x),\qquad U=\int_X^\oplus U_x\, d\mu(x),
{    \quad \mbox{ and }  \quad
  \cA=\int_X^\oplus B(\cH_x) d\mu(x).}\]
We have a measurable decomposition $X = X_0 \dot\cup X_1$, where 
$\dim \cH_x = 1$ for $x \in X_0$ and the representations
$(U_x)_{x \in X_0}$ are trivial. For $x \in X_1$, the algebras
$\cM_x$ are factors of type III$_1$ and
$(\cM_x, \Omega_x, \uline U_x)$ satisfies {\rm (Reg)}
and {\rm(Mod)}, where $\uline U_x$ is the representation of
$G/\ker(U_x)$ induced by $U_x$. Further details in \cite[Thm.~5.22]{MN24}. 

\section{An outlook on hermitian Lie algebra}

We have seen that various aspects of the geometric framework of Algebraic Quantum Field Theory (AQFT) can be explored at the Lie theory level.
In particular, the Bisognano--Wichmann and modular reflection properties give rise to a generalized wedge-boost correspondence, where the Euler elements assume a central role. As we have seen along these lines, this framework allows for the deduction of properties of wedge regions, wedge symmetries, wedge standard subspaces and wedge von Neumann algebras.

Within this generalized setting, conformal models stand out due to the significant role played by orthogonal wedges, offering a rich structure \cite{MN21}. Here the Lie algebra of the symmetry group is assumed to be simple and hermitian. Specifically a simple Lie algebra is hermitian if the center of a maximal compactly embedded subalgebra $\fk\subset\fg$ is non-zero. Equivalently, a simple real Lie algebra $\fg$ is hermitian if and only if it has a convex $\Inn(\fg)$-invariant cone $C\neq\{0\}$ or $\fg$. Furthermore a simple hermitian real Lie algebra contains an Euler element if and only if it is of tube type, i.e.~the restricted root system is of type $C_r$ (see \cite[Prop.3.11]{MN21}).  In view of the classification in Sect.\ref{sect:ee}, in the hermitian simple case there exists a unique $\Inn(\fg)$-orbit of symmetric Euler elements.  By \cite[Theorem 3.13]{MN21}, in a simple Lie algebra symmetric Euler elements have orthogonal Euler element partners: We shall say that $h_a,h_b\in \cE(\fg)$ in a Lie algebra $\fg$ are orthogonal Euler elements if they satisfy the condition $$\tau_{h_a}(h_b)=-h_b.$$ As examples, this is the case of the Euler elements associated to the right half-circle and the upper half-circle of the wedge geometry on $S^1$  and the case of Euler elements associated to $\RR^{1,d}$-Minkowksi wedges $W_i=\{x\in\RR^{1,d}, |x_0|<x_i\}$ with $i=1,2$ and $d\geq2$.

One can consider couples of orthogonal elements in a simple Hermitian Lie algebra and study orbits of orthogonal elements. The existence of disjoint orbits is particularly evident in Conformal Field Theory on Minkowksi spacetime. Given  a unitary representation $U$ of the universal covering of the conformal group~$\tilde G$, where $\fg=\so(2,d+1)$, inequivalent couples of orthogonal Euler elements  generate inequivalent representations of $\tilde \SL_2(\RR)$ not all of them satisfying a non-trivial spectral condition. A systematic analysis of these orthogonal pairs, along with a geometric examination of the conformal models, will be discussed in  forthcoming papers together with an analysis of their relationship with conformal AQFT models (\cite{MNO25, MN25}).

\end{document}